\documentclass[a4paper,10pt,twoside]{article}
\usepackage{pst-fill,pst-grad}
\usepackage{textcomp}
\usepackage[english]{babel}
\usepackage[utf8]{inputenc}
\usepackage{graphicx}
\usepackage{amsmath}
\usepackage{float}
\usepackage{fancyhdr}
\usepackage[matrix,arrow,curve]{xy}
\usepackage{pstricks} 
\usepackage{amsmath,amsfonts,verbatim,afterpage,theorem,euscript,mathrsfs,amssymb}
\usepackage{amsfonts}
\usepackage{amssymb}
\usepackage{array}
\usepackage{dsfont}
\usepackage{hyperref}

\def\vu{\vec{u}}
\def\vv{\vec{v}}
\def\vw{\vec{w}}
\def\vf{\vec{f}}
\def\vF{\vec{F}}
\def\vn{\vec{\nabla}}
\def\vor{\vec{\omega}}
\def\Div{\text{\rm div}\,}
\def\Curl{\text{\rm curl}\,}

\newtheorem{Definition}{Definition}[section]
\newtheorem{Proposition}{Proposition}[section]
\newtheorem{Lemme}{Lemma}[section]
\newtheorem{Theoreme}{Theorem}
\newtheorem{Corollaire}{Corollary}[section]
\newtheorem{Remarque}{Remark}[section]
\title{\bf The role of the pressure in the partial regularity theory for weak solutions of
the Navier--Stokes equations}
\author{Diego Chamorro\footnote{Laboratoire de Math\'ematiques et Mod\'elisation d'Evry (LaMME) - UEVE, UMR CNRS 8071 \& ENSIIE. Universit\'e d'Evry Val d'Essonne, IBGBI, 23 Boulevard de France, 91037 Evry Cedex, France}$^{\; \;, }$\footnote{email: \text{diego.chamorro@univ-evry.fr}}$\,$, Pierre-Gilles Lemari\'e-Rieusset$^{*,}$\footnote{email: \text{plemarie@univ-evry.fr}}$\, $, Kawther Mayoufi$^{*}$} 
\begin{document}
\maketitle

%%%%%%%%%%%%%%%%%%%%%%%%%%%%%%%%%%%%%%%%%%%%%%
\begin{scriptsize}
\abstract{We study the role of the pressure in the partial regularity theory for weak solutions of the Navier--Stokes equations. By introducing the notion of dissipative solutions, due to Duchon \& Robert \cite{DR1}, we will provide a generalization of the Caffarelli, Kohn and Nirenberg theory. Our approach gives a new enlightenment of the role of the pressure in this theory in connection to Serrin's local regularity criterion.}\\[3mm]
\textbf{Keywords: Navier--Stokes equations;  Partial regularity; Caffarelli, Kohn and Nirenberg theory;  Serrin criterion.}
\end{scriptsize}

%\tableofcontents
%%%%%%%%%%%%%%%%%%%%%%%%%%%%%%%%%%%%%%%%%%%%%%
%%%%%%%%%%%%%%%%%%%%%%%%%%%%%%%%%%%%%%%%%%%%%%
\section{Introduction and presentation of the results}
In this article we want to study some problems related with the role of the pressure in the \emph{partial} regularity theory for weak solutions of the Navier--Stokes equations. Before going into any further details, it is worth noting that if we work in the whole space it is possible to get rid of the (unknown) pressure in a straighforward way. Indeed, if we consider the problem to find a weak solution in $\mathbb{R}\times\mathbb{R}^3$ of the equation
\begin{equation}\label{Navier-Stokes_1}
\begin{cases}
\partial_t \vu = \nu \Delta \vu-(\vu \cdot \vn)\vu -\vn p+\vf, \quad  \Div(\vu)=0, \\[2mm]
\vu(0, x)=\vu_0 \in L^2(\mathbb{R}^3),\, \Div(\vu_0)=0,
\end{cases}
\end{equation}
where the viscosity $\nu>0$ is a fixed parameter,  the force $\vf$ is fixed in $L^2([0,+\infty[, H^{-1}(\mathbb{R}^3))$,  and  the solution $\vu$ satisfies $\vu \in L^\infty([0,T[, L^2(\mathbb{R}^3))\cap L^{2}([0,T[, \dot{H}^1(\mathbb{R}^3))$ for every $T>0$, it is then possible to use the Leray projector $\mathbb{P}$, defined by $\mathbb{P}(\varphi)=\varphi-\vn \frac{1}{\Delta}(\vn \cdot \varphi)$, to prove  that the previous problem is equivalent to the following one (see \cite{LEMA1}, Chapter 11)
\begin{equation*}
\partial_t \vu = \nu \Delta \vu-\mathbb{P}\big(\vn\cdot(\vu \otimes\vu)\big)  +\mathbb{P}(\vf),\qquad \Div(\vu)=0.
\end{equation*}
Thus, if we are interested in studying the regularity of weak solutions $\vec u$ of such problem, the pressure will not play a fundamental role. However, in a \emph{local} framework it is not possible to use this technique as the Leray projector is a non-local operator. To overcome this issue we have at our disposal two different approaches. 
%%%%%%%%%%%%%%%%%%%%%%%%%%%%%%%%%%%%%%%%%%%%%%
\subsubsection*{The Serrin regularity theory}

Following \cite{Serrin1}, \cite{Struwe} or \cite{Taka}, the first approach consists in taking the curl of equation (\ref{Navier-Stokes_1}) in order to get rid of the pressure. Then, denoting by $\vor=\Curl  \vu= \vn \wedge \vu$ and since we have $\vn \wedge \vn p\equiv 0$, we obtain the following equation:
\begin{equation}\label{Equa_Vorticity0}
\partial_{t}\vor=\nu\Delta \vor-\vn \wedge \big((\vu \cdot \vn) \vu\big)+\vn \wedge\vf,
\end{equation}
from which it is possible to study the local regularity of the vorticity $\vor$ and then to deduce the regularity of the weak solutions $\vu$. Here again, just as in the non-local case, we observe that the pressure plays no particular role.\\ 

This technique was first developped by Serrin \cite{Serrin1} in the following setting: let $Q=]a,b[\times B_{x_{0},r_{0}}$ be a bounded set where $]a,b[$ is an interval of the real line and $B_{x_{0},r_0}$ stands for the Euclidean ball $B_{x_{0},r_0}= B(x_{0},r_0)$  with $x_{0}\in \mathbb{R}^{3}$ and $r_0>0$. Let moreover $\vf\in L^2_{t}H^k_{x}(Q)$ for some $k\geq 0$, let $\vu\in L^{\infty}_{t}L^{2}_{x}(Q)\cap L^{2}_{t}\dot{H}^{1}_{x}(Q)$ and $p\in \mathcal{D}'(Q)$; if we assume that $\vu$ is a weak solution on $Q$ of the Navier--Stokes equations (\ref{Navier-Stokes_1}) then, if $\vu\in L^{\infty}_{t}L^{\infty}_{x}(Q)$, we may conclude that locally the regularity of $\vu$ is in fact driven by the regularity of $\vf$:  for every $a<c<b$ and $0<\rho<r_0$ we have that $\vu \in L^{\infty}\big(]c,b[, H^{k+1}(B_{x_{0},\rho})\big)\cap L^{2}\big(]c,b[,\dot{H}^{k+2}(B_{x_{0},\rho})\big)$.\\ 

This type of results is known as \emph{local} regularity theorems. We make now several remarks concerning the Serrin regularity criterion.
\begin{itemize}
\item[$(i)$] The pressure $p$ can be a very general object as we only need that $p\in \mathcal{D}'(Q)$ and this is not a problem since it has disappeared in equation (\ref{Equa_Vorticity0}).
\item[$(ii)$] It is clear from equation (\ref{Equa_Vorticity0}) that the regularity of $\vu$ is related to the regularity of the force $\vf$.
\item[$(iii)$] We observe that by this method it is not possible to obtain any information of the regularity in the \emph{time} variable: indeed, Serrin gave the following example: if $\phi$ is a bounded function on $\mathbb{R}$ and if $\psi$ is a harmonic function on $\mathbb{R}^{3}$, we define $\vu$ on $]0,1[\times B(0,1)$ by $\vu(t,x)=\phi(t)\vn \psi(x)$. We have $\Div\vu=\phi(t)\Delta \psi(x)=0$, $\Curl \vu=\vn \wedge \vu=0$, $\Delta \vu=0$ and we obtain $\vu\cdot \vn \vu=\vn\left(\frac{|\vu|^{2}}{2}\right)$. Now, if $\vu$ satisfies the Navier--Stokes equations (with a null force) we have $\partial_{t}\vu=\nu \Delta \vu-(\vu \cdot \vn) \vu-\vn p=-\vn\left(\frac{|\vu|^{2}}{2}\right)-\vn p$ from which we derive the following relationship $$p(t,x)=-\frac{|\vu(t,x)|^{2}}{2}-\partial_{t}\phi(t) \psi(x).$$ 
We can thus see that a control of $p$ is equivalent to a control of $\partial_{t}\vu$. Moreover we have $\vu\in L^{\infty}_{t}L^{2}_{x}\cap L^{2}_{t}\dot{H}^{1}_{x}\cap L^{\infty}_{t}L^{\infty}_{x}$ on $]0,1[\times B(0,1)$, but it is easy to see that if $\phi$ is not regular then there is no hope to obtain regularity for $\vu$ with respect to the time variable.
\item[$(iv)$] The boundedness assumption $\vu\in L^{\infty}_{t}L^{\infty}_{x}(Q)$ can be generalized: it is enough to assume that $\vec u$ has some (sub)critical behavior with respect to the scaling of the equation.  Serrin \cite{Serrin1} proved that, if $\vf \in L^2_{t}H^1_{x}(Q)$ and if $\vu\in L^p_{t} L^q_{x}(Q)$ with $\frac{2}{p}+\frac{3}{q}<1$, then,  for every $a<c<b$ and $0<\rho<r_{0}$ we  actually have that $\vu\in L^\infty_{t} L^\infty_{x}(]c,b[\times B_{x_{0}, \rho})$.  Significant efforts have been made to generalize even more this hypothesis, see for example \cite{Struwe}, \cite{Taka} or \cite{Chen}. In particular, parabolic Morrey-Campanato spaces were used by O'Leary \cite{OLeary1} to generalize Serrin's theorem and we will see how to exploit this framework later on.
\item[$(v)$] Our last remark is that Serrin's theory relies on  the subcriticality assumption $\vu\in L^p_{t}L^q_{x}(Q)$ with $\frac{2}{p}+\frac{3}{q}<1$ (or the criticality assumption  $\vu\in L^p_{t}L^q_{x}(Q)$ with $\frac{2}{p}+\frac{3}{q}=1$ and $3<q\leq +\infty$, proved by Struwe \cite{Struwe} and Takahashi \cite{Taka}), which  is indeed very restrictive: from the information that $\vu\in L^{\infty}_{t}L^{2}_{x}(Q)\cap L^{2}_{t}\dot{H}^{1}_{x}(Q)$, using the Sobolev inequalities, we can only obtain that $\vu\in L^p_{t}L^q_{x}(Q)$ with $\frac{2}{p}+\frac{3}{q}=\frac{3}{2}$. Thus, we have a supercritical behavior of $\vu$ and we cannot deduce from the usual hypotheses the Serrin criterion for local regularity.
\end{itemize}
%%%%%%%%%%%%%%%%%%%%%%%%%%%%%%%%%%%%%%%%%%%%%%
\subsubsection*{Caffarelli, Kohn and Nirenberg theory}

To circumvene this supercriticality of $\vu$, Caffarelli, Kohn and Nirenberg \cite{CKN} introduced another approach which is actually satisfied by $\vu$ in the neighborhood of almost every point of $Q=]a,b[\times B_{x_0,r_0}$, so that the lack of regularity is concentrated on a very small set. As we will work here only on neighborhoods of points, the results associated to this theory are denoted by \emph{partial} regularity theorems.\\

Let us be more precise on the Caffarelli--Kohn--Nirenberg theory and we introduce now the main ingredients of this theory
\begin{enumerate}
\item[$\bullet$] \emph{The notion of weak solutions:} a weak solution $(\vu, p)$ of the Navier--Stokes equations on a domain $Q=]a,b[\times B_{x_0,r_0}$  is a time-dependent vector field  $\vu\in L^{\infty}_{t}L^{2}_{x}(Q)\cap L^{2}_{t}\dot{H}^{1}_{x}(Q)$ and a pressure $p\in \mathcal{D}'(Q)$. If we assume that $\vu$ is a weak solution on $Q$ of the Navier--Stokes equations such that for some (unknown) pressure $p$  and some (given) force $\vec f$ we have
$$\partial_t \vu = \nu \Delta \vu-(\vu \cdot \vn)\vu -\vn p+\vf,\qquad \Div(\vu)=0,$$
then the terms  $\partial_t \vu$, $\nu \Delta \vu$ and $(\vu \cdot \vn)\vu$ are well-defined in $\mathcal{D}'(Q)$ when  $\vu\in L^{\infty}_{t}L^{2}_{x}(Q)\cap L^{2}_{t}\dot{H}^{1}_{x}(Q)$, so that the equations are meaningful for $p$ and $\vec f$ in  $\mathcal{D}'(Q)$.
\item[$\bullet$] \emph{The set of singular points:} following Serrin, we shall say that a point $(t_{0},x_{0})\in Q$ is a \emph{regular} point of the solution $\vec u$ if there exists a neighborhood $\mathcal{V}$ of $(t_{0},x_{0})$ such that $\vu\in L^\infty_{t} L^\infty_{x}(\mathcal{V})$. We then define the set $\Sigma$ of singular points as the set of points $(t,x)\in Q$ that are not regular points of $\vu$.

Observe from the remarks above that if $\vf$ is regular enough (for instance, $\vf\in L^2_{t}H^1_{x}$) it is equivalent to ask that $\vec u$ satisfies $\vec u\in L^p_tL^q_x$ on a neighborhood of $(t_{0},x_{0})$, for some $(p,q)$ with $\frac{2}{p}+\frac{3}{q}\leq 1$ and  $3<q\leq +\infty$.  
\item[$\bullet$] \emph{The set of large gradients for the velocity:} we shall say that a point $(t,x)\in Q$ is a point of large gradients for the velocity if we have
\begin{equation*}
\underset{r\to 0^{+}}{\lim \sup}\frac{1}{r}\iint_{]t-r^{2}, t+r^{2}[\times B_{x,r}}|\vn \otimes \vu|^{2}ds\, dy>0.
\end{equation*}
We shall write $\Sigma_0$ for the set of points of large gradients.\\

Remark that if $\vf$ is regular enough (for instance, if $\vf\in L^2_{t}L^2_{x}(Q)$), then Serrin's analysis tells us that, for a regular point $(t_{0},x_{0})\notin \Sigma$, there exists a neighborhood $\mathcal{V}$ of this point such that $\vu\in L^\infty_t H^1_x$ on $\mathcal{V}$: thus, for $r$ small enough we have 
$$\iint_{]t_{0}-r^{2}, t_{0}+r^{2}[\times B_{x_{0},r}}|\vn \otimes \vu|^{2}ds\, dy=O(r^2),\quad \mbox{and then }(t_{0},x_{0})\notin \Sigma_0.$$ 
Hence, $\Sigma_0$ is in fact a set of singular points. Besides, $\Sigma_0$ is a very small set: indeed, Caffarelli, Kohn and Nirenberg \cite{CKN}  showed how to deduce from the hypothesis $\vu\in L^2_t H^1_x$ the fact that the (parabolic)  one-dimensional Hausdorff measure $\mathcal{H}^{1}_{2}$ of $\Sigma_0$ is null. The aim of a partial regularity theory is then to find criteria that ensure that $\Sigma=\Sigma_0$, so that there are very few singular point (if any).
 \item [$\bullet$]\emph{The notion of suitable solutions:} the key point in partial regularity theory for Navier--Stokes equations is the local energy inequality first studied by Scheffer \cite{Scheffer1, Scheffer2}. If $\vec f$ and $p$ are regular enough to ensure that the products $p\vec u$ and $\vec f\cdot\vec u$ are meaningful as distributions, then the quantity
\begin{equation}\label{Defmu} 
\mu=-\partial_t|\vu|^{2}+\nu\Delta|u|^{2}-2\nu|\vec\nabla\otimes\vec u|^2-\Div(|\vu|^{2}\vu)-2\Div(p \vu)+2\vf \cdot \vu,
\end{equation} 
is well-defined as a distribution in $\mathcal{D}'(Q)$. 
\begin{Remarque}\label{Remarque0}
For the force Kukavica \cite{K} proposed that $\vec f\in L^{10/7}_{t}L^{10/7}_{x}(Q)$, as $\vec u\in L^{10/3}_{t}L^{10/3}_{x}(Q)$ due to the Sobolev embedding inequalities. For the pressure $p$ Vasseur \cite{Vasseur} showed that $p\in L^r_t L^1_x(Q)$ with $r>1$ was enough.
\end{Remarque}

Moreover, if $\vec u$ is regular enough, we will define {\em suitable solutions} as the solutions for which inequality 
$$\mu\geq 0,$$
holds, \emph{i.e.} for which the distribution $\mu$ is given by a locally finite non-negative Borel measure.\\ 

Remark in particular that if we know $\vec u\in L^p_{t}L^q_{t}(Q)$ with $p=q=4$ (which is not in the scope of the Serrin criterion since in this case we have $\frac{2}{p}+\frac{3}{q}>1$) then we have $\mu=0$: we are indeed in a more general framework. See \cite{LEMA2} for a proof of this fact.\\
\end{enumerate}
Once we have detailed the setting, we can state the  Caffarelli--Kohn--Nirenberg regularity theorem: let $\Sigma_\epsilon$ be the set of points $(t,x)$ in $Q$ such that
\begin{equation*}\label{Condition_epsilon}
\underset{r\to 0^{+}}{\lim \sup}\frac{1}{r}\iint_{]t-r^{2}, t+r^{2}[\times B_{x,r}}|\vn \otimes \vu|^{2}ds\, dy>\epsilon.  
\end{equation*} 
Then Caffarelli, Kohn and Nirenberg \cite{CKN} proved that $\Sigma=\Sigma_0=\Sigma_\epsilon$ for some $\epsilon>0$ small enough that does not depend on $Q$, $\vf$ nor $\vu$ provided the following assumptions are fulfilled :
\begin{itemize}
\item[$\ast$] $p$ is regular enough (usually, one takes $p\in L^{3/2}_{t}L^{3/2}_{x}(Q)$ see \cite{Lin}, however  Vasseur \cite{Vasseur} showed that $p\in L^r_t L^1_x(Q)$ with $r>1$ was enough),
\item[$\ast$] $\vf$ is regular enough (in \cite{CKN}, the condition was $\vf \in L^\rho_{t}L^\rho_{x}(Q)$ with $\rho>5/2$, but  other assumptions can be made, in particular $\vf\in L^2_t H^1_x(Q)$ will be enough),
\item[$\ast$] $\vec u$ is suitable, \emph{i.e.}, the associated distribution $\mu$ defined by identity (\ref{Defmu}) is non-negative.\\
\end{itemize}

As a matter of fact, Caffarelli, Kohn and Nirenberg proved a slightly more general result: under  regularity assumptions on $p$ (such as $p\in L^{3/2}_{t,x}$) and on $\vf$ (such as $\vf\in L^{5/2+\varepsilon}_{t,x}$), and under the suitablity assumption on $\vu$ (i.e. $\mu\geq 0$), then the solution $\vu$ is H\"older-regular on \emph{both} time and space variables in a neighborhood $\mathcal{V}$ of $(t,x)$: for some $\eta>0$ and $C\geq 0$, we have, for $(s,y)\in \mathcal{V}$ and $(\tau,z)\in \mathcal{V}$, 
$$\vert \vec u(s,y)-\vu(\tau,z)\vert\leq C (\sqrt{\vert s-\tau\vert}+\vert y-z\vert)^\eta.$$

Several remarks are in order here.
\begin{itemize}
\item[$(i)$]  First, we notice that it is necessary to impose some conditions in the pressure $p$ and in this sense this approach is less general than the Serrin criterion where we only have $p\in \mathcal{D}'$. 
\item[$(ii)$]  We do not need the local boundedness assumption $\vu\in (L^{\infty}_{t}L^{\infty}_{x})_{loc}$ (or $(L^{p}_{t}L^{q}_{x})_{loc}$ with $\frac{2}{p}+\frac{3}{q}\leq 1$ and $3<q\leq +\infty$); instead we will use the hypothesis of \emph{suitability}. In this sense the Caffarelli-Kohn-Nirenberg theory is more general, at least when we are studying the constraints on $\vu$. 
\item[$(iii)$] Since we have some control in the pressure $p$ it is possible to obtain some regularity in the time variable as it was underlined by the Serrin example.
\item[$(iv)$] Finally, it is worth noting that this regularity is only obtained in small neighborhoods of points.
\end{itemize}
As we can see, if we compare the hypotheses and the conclusions of these theories that study the local/partial regularity for the weak solutions of the Navier--Stokes equations we obtain two quite different approaches.
%%%%%%%%%%%%%%%%%%%%%%%%%%%%%%%%%%%%%%%%%%%%%%
\subsubsection*{Presentation of the results}

The general aim of this article is to weaken as much as possible the regularity assumption on the pressure, while keeping the main lines of the Caffarelli--Kohn--Nirenberg theory and, by doing so, we will obtain a different point of view to the partial regularity theory. \\

We use as a starting point the following remark: the Caffarelli--Kohn--Nirenberg theory is based on the suitability of the solution $\vu$, \emph{i.e.} on the  local energy inequality $\mu\geq 0$, where $\mu$ is given by identity (\ref{Defmu}), and we can see that indeed we have two assumptions in the definition of \emph{suitability}:
\begin{itemize}
\item the pressure $p$ is regular enough to allow the quantity $\mu$ to be defined as a distribution,
\item and the fact that $\mu$ is non-negative.
\end{itemize}
But if we just assume $p\in \mathcal{D}'(Q)$, we have that $\mu$ is no longer well-defined and thus we must change the definition of suitability. So, our first task is to give a sense to the product $p\vu$ even when $p\in \mathcal{D}'(Q)$:
%%%%%%%%%%%%%%%%%%%%%%%%%%%%%%%%%%%%%%%%%%%%%%
\begin{Proposition}\label{Proposition_40} Let $x_{0}\in \mathbb{R}^{3}$ and $\rho>0$, we consider $Q=]a,b[\times B_{x_{0},\rho}$ a bounded subset of $\mathbb{R}\times \mathbb{R}^{3}$. Assume that $\vu\in L^{\infty}_{t}L^{2}_{x}(Q)\cap L_{t}^{2}\dot{H}_{x}^{1}(Q)$ with $\Div(\vu)=0$ and $p\in \mathcal{D}'(Q)$ are solutions of the Navier--Stokes equations on  $Q$ :
$$\partial_{t} \vu=\nu\Delta \vu -(\vu \cdot \vn )\vu -\vn p+\vf,$$
where $\vf\in L^{10/7}_t L^{10/7}_x(Q)$ and $\Div(\vf)=0$. \\

Let $\gamma\in \mathcal{D}(\mathbb{R})$ and  $\theta \in \mathcal{D}(\mathbb{R}^{3})$ be two smooth functions such that $\displaystyle{\int_{\mathbb{R}}} \gamma(t)dt=\displaystyle{\int_{\mathbb{R}^{3}}} \theta(x)dx=1$, $supp(\gamma)\subset ]-1,1[$ and $supp(\theta)\subset B(0,1)$. We set for $\alpha, \varepsilon >0$ the functions $\gamma_{\alpha}(t)=\frac{1}{\alpha}\gamma(\frac{t}{\alpha})$ and $\theta_{\varepsilon}=\frac{1}{\varepsilon^{3}}\theta(\frac{x}{\varepsilon  })$ and we define 
$$\varphi_{\alpha, \varepsilon}(t,x)=\gamma_{\alpha}(t)\theta_{\varepsilon}(x).$$ 
Then, if the set $Q_{t_0,x_0,r_{0}}=]t_{0}-r_{0}^{2}, t_{0}+r_{0}^{2}[\times B_{x_{0}, r_{0}}$ is contained in $Q$, the distributions $\vu \ast \varphi_{\alpha, \varepsilon}$ and\footnote{Convolutions are considered in the time \emph{and} the space variable.} $p\ast \varphi_{\alpha, \varepsilon}$ are well defined in the set $Q_{t_0,x_0,r_{0}/4}\subset Q$ for $0<\alpha<r_0^2/2 $ and $0<\varepsilon<r_0/2$. Moreover,  the limit
$$\underset{\varepsilon \to 0}{\lim}\,\underset{\alpha \to 0}{\lim}\, \Div\left[\big(p\ast \varphi_{\alpha, \varepsilon}\big) \times \big(\vu \ast \varphi_{\alpha, \varepsilon}\big)\right],$$
exists in $\mathcal{D}'$ and does not depend on the choice of $\gamma$ and $\theta$. 
\end{Proposition}
%%%%%%%%%%%%%%%%%%%%%%%%%%%%%%%%%%%%%%%%%%%%%%
We shall write from now on
\begin{equation*}
\underset{\varepsilon \to 0}{\lim}\,\underset{\alpha \to 0}{\lim}\, \Div\left[\big(p\ast \varphi_{\alpha, \varepsilon}\big) \times \big(\vu \ast \varphi_{\alpha, \varepsilon}\big)\right]= \langle\Div(p\vu)\rangle.
\end{equation*}
The existence of this limit is not absolutely trivial but it will allows us to work with the object $\langle\Div(p\vu)\rangle$ where the pressure $p$ belongs to $\mathcal{D}'(Q)$: we can now introduce the following concept that will replace the notion of suitability.
%%%%%%%%%%%%%%%%%%%%%%%%%%%%%%%%%%%%%%%%%%%%%%
\begin{Definition}[Dissipative solutions]\label{Definition_Dissipative} Within the framework of the Proposition \ref{Proposition_40}, i.e.:
\begin{itemize}
\item assume $\vf\in L^{10/7}_t L^{10/7}_x(Q)$ with $Q=]a,b[\times B_{x_{0},\rho}$ where $x_{0}\in\mathbb{R}^{3}$, $\rho>0$ and $\Div(\vf)=0$,
\item  if $(\vu,p)$ is a solution of the Navier--Stokes equations $\partial_{t} \vu=\nu\Delta \vu -(\vu \cdot \vn )\vu -\vn p+\vf$, $\Div(\vu)=0$ on $Q$ with $\vu\in L^{\infty}_{t}L^{2}_{x}(Q)\cap L_{t}^{2}\dot{H}_{x}^{1}(Q)$ and $p\in \mathcal{D}'(Q)$,
\end{itemize} 
we will say that a solution $\vu$ is  \emph{dissipative} if the distribution $M$ given by the expression
\begin{equation}\label{Dissipative1}
M=-\partial_{t}|\vu|^{2}+\nu \Delta |\vu|^{2}-2\nu|\vn \otimes \vu|^{2}-\Div(|\vu|^{2}\vu)-2\langle\Div(p\vu)\rangle+2 \vf\cdot \vu,
\end{equation}
is a non-negative locally finite measure on $Q$.
\end{Definition}
%%%%%%%%%%%%%%%%%%%%%%%%%%%%%%%%%%%%%%%%%%%%%%
We have here two remarks.
\begin{itemize}
\item[$(i)$] It is clear since $p\in \mathcal{D}'$ that the concept of dissipative solutions is more general than notion of suitable solutions, and we will show that suitability implies dissipativity. Moreover, it is possible to show that the set of dissipative solutions is strictly larger than suitable solutions one, see Remark \ref{Remarque1} below for the details.
\item[$(ii)$] It is worth noting that if we assume $\vu \in L^{4}_{t}L^{4}_{x}$ (with $\vf$ regular enough, say $\vf\in L^{10/7}_t L^{10/7}_x$, but without any condition on $p$) then we can prove that  we actually have $M=0$. See Remark \ref{Remarque44} below. \\
\end{itemize}

Once we have introduced the notion of dissipative solutions, we now may state our main theorem:
%%%%%%%%%%%%%%%%%%%%%%%%%%%%%%%%%%%%%%%%%%%%%%
\begin{Theoreme}\label{Theorem_1}
Let $\nu>0$ be a fixed parameter. Let $Q=]a,b[\times B_{x,\rho}$ be a bounded domain of $\mathbb{R}\times \mathbb{R}^3$ and let $(\vu, p)$ be a weak solution on $Q$ of the Navier--Stokes equations
\begin{equation}\label{Equation_NS}
\begin{cases}
\partial_t \vu = \nu \Delta \vu-(\vu \cdot \vn)\vu -\vn p+ \vf, \quad div(\vu)=0,\\[2mm]
\vu(0, x)=\vu_0 \in L^2(\mathbb{R}^3), \, div(\vu_0)=0.
\end{cases}
\end{equation}
We assume that: 
\begin{itemize}
\item[$\bullet$] $\vu \in L^{\infty}_t L^2_x(Q) \cap L^2_t \dot{H}^1_x (Q)$ and $p\in \mathcal{D}'(Q)$,\\ 
\item[$\bullet$] $\vf \in L^{2}_{t}H^{1}_{x}(Q)$,\\
\item[$\bullet$] $\vu$ is dissipative in the sense of the Definition \ref{Definition_Dissipative} given above.\\
\end{itemize}
There exists a positive constant $\varepsilon^{*}>0$, which depends only on $\nu$, such that, if for some point $(t_0, x_0)\in Q$ we have the inequality
\begin{equation}\label{Smallness_Condition_1}
\underset{r\to 0}{\limsup} \frac{1}{r}\underset{]t_0-r^2, t_0+r^2[\times B_{x_0, r}}{\int \int} |\vn \otimes \vu (t,x)|^2 dt\, dx<\varepsilon^{*},
\end{equation}
then the solution $\vu$ is bounded in a neighborhood of $(t_0, x_0)$. In particular the point $(t_{0}, x_{0})$ is regular. 
\end{Theoreme}
%%%%%%%%%%%%%%%%%%%%%%%%%%%%%%%%%%%%%%%%%%%%%%
It is interesting to contrast our approach to the Serrin and to the Caffarelli--Kohn--Nirenberg theories: from the point of view of the hypotheses on $\vu$ we only assume the smallness condition (\ref{Smallness_Condition_1}) and we require the dissipativeness property. 
Thus, since we impose less conditions over the pressure and since the concept of dissipative condition is more general than the suitable one, we could think our method as a generalization of the Caffarelli--Kohn--Nirenberg theory. But, as we only require that $p\in \mathcal{D'}$, due to the Serrin example given at the beginning of this article, it is not possible to expect any regularity in the time variable and the conclusion of our main theorem only provides regularity in the space variable: in this sense our result should also be considered as a generalization of the Serrin theory since we require less conditions over $\vu$.\\ 

As we can see, our method is a generalization of these theories following very specific directions: a weaker control over the pressure generates the loss of the time regularity but it is still possible to obtain regularity in the space variable.\\

Let us explain now the global strategy of the proof that will be displayed for Theorem \ref{Theorem_1}. First, as we do not make any particular assumption over the pressure $p$ (recall that we only have $p\in \mathcal{D}'$) we will take the curl in the Navier--Stokes equations (\ref{Equation_NS}) and doing so we will immediately get rid of the pressure $p$. However we will not going to work with the variable $\vor=\vn \wedge \vu$ and the corresponding equation (\ref{Equa_Vorticity0}), in fact we want to work with a more regular distribution. Furthermore, as we are interested to study the regularity problem only in a neighborhood of a point $(t_{0}, x_{0})$, we will first introduce a cutting function $\psi \in \mathcal{D}(\mathbb{R}\times \mathbb{R}^{3})$ which is null outside a small ball centered in the point $(t_{0}, x_{0})$ and then we will define a new variable in the following way:
$$\vv =  -\frac{1}{\Delta} \vn \wedge (\psi \vn \wedge \vu).$$
The crucial point here is that, roughly speaking, the two successive derivatives given by the gradients $\vn$ in the previous formula are \emph{locally compensated}, in a sense that will be made precise later on, by the operator $ -\frac{1}{\Delta}$ and thus, some of the properties of $\vu$ will be very similar to those of the new variable $\vv$ and viceversa: actually we will see that (locally) $\vec v$ is equal to $\vu$ up to a harmonic (in the space variable) correction. \\

The central idea is then to use $\vv$ as a support function to study the regularity properties of $\vu$ and thus our first task will be to describe some of the properties of $\vv$.\\

In Proposition \ref{Proposition_1} we give some basic properties of $\vv$ that can be easily deduced from the hypotheses on $\vu$; however, the full strength of this new variable will appear clearly when we will study the equations satisfied by $\vv$: indeed, we will see that the function $\vv$ satisfies the following Navier--Stokes equation (called the \emph{companion equation})
$$\partial_t \vv =\Delta \vv -  (\vv \cdot\vn) \vv- \vn q+\vF,$$
where the pressure $q$ and the force $\vF$ will be deduced from the original parameters and we will see that $q$ and $\vF$ satisfy some interesting properties. It is worth noting here that $q$ and $\vF$ will not depend on the pressure $p$. See Proposition \ref{Proposition_2} below.\\

With the help of the variable $\vv$ and the companion equation we will prove Proposition \ref{Proposition_40} and we will see that if $\vu$ is a dissipative solution then our new variable $\vv$ is actually \emph{suitable} is the sense of the Caffarelli--Kohn--Nirenberg theory. Moreover, we will prove that if $\vu$ satisfies the smallness assumption (\ref{Smallness_Condition_1}) then it will be also the case for the new variable $\vv$. Thus, since the function $\vv$ satisfies the Navier--Stokes equations in a slightly different framework than the function $\vu$, it would be possible to study further local properties of the variable $\vv$. Then we will see how these properties of the function $\vv$ are transmitted to the original variable $\vu$. Finally, the last step which is given with the Proposition \ref{Proposition_403} will explain how to deduce that $\vu$ is a locally bounded function and Theorem \ref{Theorem_1} will be completely proven.\\

The plan of the article is the following. In Section \ref{Secc_Parabolic} we recall some of the tools that will be used in the proof of Theorem \ref{Theorem_1}, in particular we will insist in the parabolic setting of the problem and in O'Leary's and Kukavica's theorems which are crucial in our study. 
Section \ref{Secc_Preuve_Theo1} is devoted to the proof of the main theorem. Technical lemmas are postponed to the appendix.
%%%%%%%%%%%%%%%%%%%%%%%%%%%%%%%%%%%%%%%%%%%%%%
\section{Parabolic scaling and related tools}\label{Secc_Parabolic}
The main idea in the proof of Caffarelli, Kohn and Nirenberg \cite{CKN} is to estimate the regularity of $\vu$ by estimating the size of some scaled integrals. Of course, it uses the invariance of the Navier--Stokes equations under a well-defined rescaling :  if $\vu$ is a solution of 
$$\partial_t \vu = \nu \Delta \vu-(\vu \cdot \vn)\vu -\vn p+\vf,\qquad \Div\vu=0,$$ 
on $]t_0-a, t_0+a[\times B_{x_0,r_0}$, then for $\lambda>0$ we have that 
$$\vec u_\lambda(t,x)=\lambda \vec u\big(t_0+\lambda^2(t-t_0), x_0+\lambda (x-x_0)\big),$$ 
is still a solution of the Navier--Stokes equations on the rescaled domain $]t_0-\frac{a}{\lambda^2},t_0+\frac{a}{\lambda^2}[\times B_{x_0,\frac{r_0}{\lambda}}$, for the rescaled force $\vec f_\lambda(t,x)=\lambda^3 \vec f(t_0+\lambda^2(t-t_0), x_0+\lambda (x-x_0))$ and the rescaled pressure $p_\lambda(t,x)=\lambda^{2}p(t_0+\lambda^2(t-t_0), x_0+\lambda (x-x_0))$.\\ 

Thus it is natural to work within the frame of the geometry generated through those parabolic scalings. Hence, we shall consider the parabolic distance on   $\mathbb{R}\times \mathbb{R}^{3}$  defined by
\begin{equation}\label{Definition_Distance}
d_{2}\big((t,x), (s,y)\big)=\max\{|t-s|^{\frac{1}{2}},|x-y|\},
\end{equation}
and we will denote by $Q_{t, x, r}$ the parabolic ball of center $(t,x)$ and radius $r$ \emph{i.e.} 
\begin{equation}\label{Def_BouleParabolique}
Q_{t, x, r}=\big\{(s,y)\in\mathbb{R}\times \mathbb{R}^{3}:  d_{2}\big((t,x), (s,y)\big)<r\big\}.
\end{equation}
The space $\mathbb{R}\times\mathbb{R}^3$, endowed  with the parabolic distance $d_2$ and the Lebesgue measure $dt\, dx$ is a space of homogeneous type (in the sense of \cite{CW}) with homogenous dimension equal to $5$:  $\displaystyle{\int_{Q_{t,x,r}}\, dt\, dx =C r^5}$. Associated to this distance, we have the notion of Hausdorff measures. If $A$ is a subset of $\mathbb{R}\times \mathbb{R}^{3}$, we define for $\delta>0$ the set $\mathcal{I}_\delta(A)$ as the set of countable families $(Q_n)_{n\in\mathbb{N}}$ of parabolic balls   $Q_n=Q_{t_n,x_n,r_n}$ such that $A\subset\underset{{n\in\mathbb{N}}} \bigcup Q_n$ and $\underset{n\in\mathbb{N}}{\sup}\; r_n<\delta$. For $\alpha>0$, we define for all $\delta>0$ the quantity
 $$\mathcal{H}_{2,\delta}^{\alpha}(A)=\inf_{(Q_n)_{n\in\mathbb{N}}\in \mathcal{I}_\delta(A)}  \left\{\sum_{n=0}^{+\infty}r_n^{\alpha} \right\}.$$
  Now, in order to obtain the Hausdorff measure $\mathcal{H}_{2}^{\alpha}$ of the set $A$ we make $\delta\longrightarrow 0$:
\begin{equation*}
\mathcal{H}_{2}^{\alpha}(A)=\underset{\delta \to 0}{\lim }\mathcal{H}_{2,\delta}^{\alpha}(A).
\end{equation*}
For more details concerning the properties of the Hausdorff measure see \cite{Federer}. The result of 
Caffarelli, Kohn and Nirenberg \cite{CKN} states precisely that the set $\Sigma_0$ of points with large gradients for the velocity satisfies $$\mathcal{H}_{2}^{1}(\Sigma_0) =0.$$

While the Hausdorff measures involves only the parabolic distance, we shall also work with   parabolic Morrey spaces which involve both the distance and the measure. Let $\mathcal{Q}$ be the collection of parabolic balls $Q_{t,x,r}$ of type (\ref{Def_BouleParabolique}) where $t\in\mathbb{R}$, $x\in\mathbb{R}^3$ and $r>0$. For  $1<p,q<+\infty$, the parabolic Morrey space $\mathcal{M}^{p,q}_2(\mathbb{R}\times \mathbb{R}^3)$ will be defined as the space 
of locally integrable functions $ {f}$ on $\mathbb{R}\times \mathbb{R}^{3}$ such that  
$$\| {f}\|_{\mathcal{M}^{p,q}_2}=\underset{ Q_{t,x,r}\in\mathcal{Q}}{\sup} \left(\frac{1}{r^{5(1-\frac{p}{q})}}\iint_{ Q_{t,x,r}}| {f}(s,y)|^{p} dsdy\right)^{\frac{1}{p}}<+\infty.$$
 Remark that $L^q(\mathbb{R}\times\mathbb{R}^3)= \mathcal{M}^{q,q}_2(\mathbb{R}\times\mathbb{R}^3)$. See the book \cite{LEMA2} for many interesting examples of applications of these functional spaces to the study of the Navier--Stokes equations.\\

Though we shall not use them in the following, it is worth recalling the related notion of parabolic  Morrey--Campanato spaces \cite{Camp}. The space  $\mathcal{L}^{p,\lambda}_2$ is defined as the space of locally integrable functions $ {f}$ on $\mathbb{R}\times \mathbb{R}^{3}$ such that  
$$\| {f}\|_{\mathcal{L}^{p,\lambda}_2}=\underset{ Q_{t,x,r}\in\mathcal{Q}}{\sup}\quad \frac{1}{r^{\lambda}} \left(\iint_{ Q_{t,x,r}}| {f}(s,y)-m_{Q_{t,x,r}}(f)|^{p} dsdy\right)^{\frac{1}{p}}<+\infty,$$ 
where $m_{Q_{t,x,r}}(f)=\displaystyle{\frac{1}{\vert Q_{t,x,r}\vert} \iint_{Q_{t,x,r}} f(s,y)\, ds\, dy}$.\\

For $0<\lambda< 5/p$, we find that $f\in \mathcal{L}^{p,\lambda}_2$ if and only if $f=g+ C$, where $C$ is a constant and $g\in \mathcal{M}^{p,q}_2$ with $\lambda= 5(\frac{1}{p}-\frac{1}{q})$. In the particular case when $5/p<\lambda<5/p+1$, we find that $f\in \mathcal{L}^{p,\lambda}_2$ if and only $f$ is H\"olderian (with respect to the parabolic distance) with H\"older regularity exponent $\eta$, where $\lambda=\frac{5}{p}+\eta$: 
$$\sup_{(t,x)\neq (s,y)} \frac{\vert f(t,x)-f(s,y)\vert}{d_2((t,x),(s,y))^\eta}<+\infty,$$
where the distance $d_{2}$ is given by formula (\ref{Definition_Distance}).\\

Morrey spaces are useful in the theory of regularity for solutions of the heat equation, hence of the Navier--Stokes equations (see for example the books \cite{LEMA1, LEMA2}). We will say that a function $f$ belongs to those spaces in a neighborhood of a point $(t_0,x_0)$ if there is a smooth compactly supported function $\varphi $ equal  to $1$ on a neighborhood of $(t_0,x_0)$ such that $\varphi f$ belongs to the Morrey space. Then a solution $h$ of the heat equation $$\partial_t h-\nu\Delta h=g+\Div(\vec H),$$
 will be locally H\"olderian of exponent $\eta$ on a neighborhood of $(t_0,x_0)$ if the data  $g$ and $\vec H$  are regular enough: on a neighborhood of $(t_0,x_0)$, we may ask that $g$ is locally $\mathcal{M}^{p_0,q_0}_2$ with $1<p_0$ and $5>q_0>5/2$ (with $\eta=2-\frac{5}{q_0}$) and $\vec H$ is locally $\mathcal{M}^{p_1,q_1}_2$   with $1<p_1$, $q_1>5$ (and $\eta=1-\frac{5}{q_1}$). 
\begin{Remarque}\label{RemarqueMorrey}
An useful remark is the following one: if $g$ has a bounded support  and belongs to $\mathcal{M}^{p,q}_2$ with $1<p\leq q<+\infty$, then $g$ belongs to $\mathcal{M}^{p_1,q_1}_2$ whenever $1<p_1\leq p$ and $p_1\leq q_1\leq q$.
\end{Remarque}
Motivated by this framework, several authors have replaced the traditional $L_{t}^pL_{x}^q$ criterion for the Serrin  local regularity result or  the Caffarelli--Kohn--Nirenberg partial regularity result by assumptions on the (local) Morrey norms of $\vec u$ and $\vec f$. In particular, O'Leary \cite{OLeary1} have stated the following variant of Serrin's regularity result which will be useful for our purposes:
%%%%%%%%%%%%%%%%%%%%%%%%%%%%%%%%%%%%%%%%%%%%%% 
 \begin{Theoreme}\label{Theoreme_2} Let $\Omega$ be a bounded domain of $\mathbb{R}^{3}$ of the form $\Omega=Q_{t_{0},x_{0},r_{0}}$ for some $t_{0}\in\mathbb{R}$, $x_{0}\in\mathbb{R}^3$ and $r_{0}>0$. Let $\vu$ be a weak solution for the Navier--Stokes equations (\ref{Equation_NS}) such that $\vu \in L^{\infty}_t L^2_x(\Omega) \cap L^2_t \dot{H}^1_x (\Omega)$, $p\in \mathcal{D}'(\Omega)$, and $\vf \in L^{2}_{t}H^{1}_{x}(\Omega)$. If moreover we have $\mathds{1}_{\mathcal{V}}\vec{u}\in\mathcal{M}^{3,\tau}_{2}(\mathbb{R}\times\mathbb{R}^{3})$, for some $\tau>5$, where $\mathcal{V}$ is a neighborhood of the point $(t_{0},x_{0})$,  then $\vu$ is a locally bounded function: for every parabolic ball $Q$ which is compactly supported in $\mathcal{V}$, we have $\vec u\in L^{\infty}_{t}L^{\infty}_{x}({Q})$.
\end{Theoreme}
%%%%%%%%%%%%%%%%%%%%%%%%%%%%%%%%%%%%%%%%%%%%%%
O'Leary stated his theorem with a null force $\vf$; however, it is not difficult to extend it to the case of a regular force $\vf\in L^2_tH^1_x$ (see \cite{LEMA2} for a proof).\\
 
Let us recall now that Caffarelli, Kohn and Nirenberg \cite{CKN} stated their theorem with regular force ($\vf\in L^\rho_{t}L^\rho_{x}$ with $\rho>5/2$) and pressure ($p\in L^{3/2}_{t}L^{3/2}_{x}$, see for instance \cite{Lin}).  Ladyzhenskaya and Seregin \cite{Lady1} then proved the Caffarelli--Kohn--Nirenberg theorem in the setting of parabolic Morrey spaces: they assumed that, on a neighborhood $\mathcal{V}$ of $(t_0,x_0)$, the force $\vf$ satisfied $\mathds{1}_{\mathcal{V}} \vf\in   \mathcal{M}^{2,q}_2$  with  $q>5/2$.\\

More recently Kukavica \cite{K}  considered less regular forces and split the proof in three steps. 
%%%%%%%%%%%%%%%%%%%%%%%%%%%%%%%%%%%%%%%%%%%%%%
\begin{Theoreme}\label{theo_kukavica} Let $\Omega$ be a bounded domain of $\mathbb{R}^{3}$ of the form $\Omega=Q_{t_{0},x_{0},r_{0}}$ for some $t_{0}\in\mathbb{R}$, $x_{0}\in\mathbb{R}^3$ and $r_{0}>0$. Let $\vu$ be a weak solution for the Navier--Stokes equations (\ref{Equation_NS}) such that $\vu \in L^{\infty}_t L^2_x(\Omega) \cap L^2_t \dot{H}^1_x (\Omega)$, $p\in \mathcal{D}'(\Omega)$, and $\vf \in  \mathcal{D}'(\Omega)$ with $\Div(\vf)=0$.  Then:
\begin{itemize}
\item[{\bf 1)}] {\bf Energy inequality.} If $p$ is regular enough (i.e., $p\in L^{q_0}_{t}L^{q_0}_{x}(\Omega)$ for some $q_0>1$) and $\vf$ is regular enough (i.e., $\vf\in L^{10/7}_{t}L^{10/7}_{x}(\Omega)$), then the quantity
 \begin{equation}\label{Def_Suitable}
\mu=-\partial_{t}|\vu|^{2}+\nu \Delta |\vu|^{2}-2\nu|\vn \otimes \vu|^{2}-\Div((|\vu|^{2}+2p)\vu)+\vec f\cdot \vec u,
\end{equation}
 is well-defined as a distribution. The solution $\vec u$ is called {\rm suitable} if $\mu$ is a locally finite non-negative measure on $\Omega$: for all $\varphi\in \mathcal{D}'(\Omega)$ such that $\varphi\geq 0$, we have
 $$ \int_{\mathbb{R}}\int_{\mathbb{R}^{3}} |\vu|^{2} (\partial_t\varphi+\nu\Delta\varphi) +(\vu \cdot \vf - 2\nu|\vn \otimes \vu|^{2}) \varphi +(|\vu|^{2}+2p)\vu\cdot\vec\nabla\varphi\, dt\, dx\geq 0.$$
\item[{\bf 2)}] {\bf The small gradients criterion.}  Assume that:
\begin{itemize}
\item[$\bullet$] $p\in L^{q_0}_{t}L^{q_0}_{x}(\Omega)$ for some $q_0>1$,
\item[$\bullet$] $\mathds{1}_{\Omega} \vf \in\mathcal{M}^{\frac{10}{7},\tau_0}_{2}$ for some $\tau_0>5/3,$
\item[$\bullet$] $\vu$ is suitable.
\end{itemize}
 There exists   positive constants $\epsilon^*>0$ and $\tau_1>5$  which depend  only on $\nu$,  $q_0$, and $\tau_0$    such that, if $(t_0,x_0)\in \Omega$ and 
 $$\limsup_{r\rightarrow 0} \frac{1}{r}\iint_{]t_0-r^2,t_0+r^2[\times B(x_0,r)} \vert\vec\nabla\otimes \vec u(s,y)\vert^2\, ds\, dy<\epsilon^*,$$ 
 then  there exists a  small parabolic neighborhood $\mathcal{Q}=Q_{t_{0}, x_{0}, \bar{r}}$ of $(t_0,x_0)$  such that we have $\mathds{1}_{\mathcal{Q}} \,\vec u\in \mathcal{M}^{3,\tau_1}_2$ and $\mathds{1}_{\mathcal{Q}}\,  p\in \mathcal{M}^{q_0,\tau_1/2}_2$. 
 
\item[ {\bf 3)}] {\bf  Regular points.} Assume that there exists a    neighborhood $\mathcal{Q}=Q_{t_{0}, x_{0}, \bar{r}}$  of $(t_0,x_0)\in \Omega$  such that 
\begin{itemize}
\item[$\bullet$]  $\mathds{1}_{\mathcal{Q}}\, \vec u\in \mathcal{M}^{3,\tau_1}_2$ for some $\tau_1>5$,
\item[$\bullet$]  $\mathds{1}_{\mathcal{Q}} \, p\in \mathcal{M}^{q_0,\tau_2}_2$ for some $1<q_0\leq \tau_2$ and $\tau_2>5/2$,
\item[$\bullet$]  $\mathds{1}_{\mathcal{Q}} \, \vf \in\mathcal{M}^{\frac{10}{7},\tau_3}_{2}$ for some $\tau_3>5/2.$  
\end{itemize}
Then there exist $0<\rho<\bar{r}$ and $\eta\in ]0,1[$ such that $\vu$ is H\"olderian (with parabolic H\"older regularity exponent $\eta>0$) on $Q_{t_0,x_0,\rho}$. In particular,  the point $(t_0,x_0)$ is regular.
\end{itemize}
\end{Theoreme}
%%%%%%%%%%%%%%%%%%%%%%%%%%%%%%%%%%%%%%%%%%%%%%
We end this section with some remarks concerning the hypotheses stated for the force $\vf$: 
\begin{itemize}
\item[$(i)$] in the first point of this theorem, we are only interested to give a sense to the product $\vf\cdot \vu$, thus since we have $\vu\in L^{10/3}_{t}L^{10/3}_{x}$ it is enough to assume that $\vf\in L^{10/7}_{t}L^{10/7}_{x}$.  
\item[$(ii)$] However, for the second point we will need more regularity and if we want to work with more classical spaces we may ask that $\vf\in L^2_t L^2_x$. Indeed, since $L^2_tL^2_x=\mathcal{M}^{2,2}_2$, and since $\mathcal{Q}$ is a bounded subset, we find that $\mathds{1}_{\mathcal{Q}}\vf\in L^2_t L^2_x$ implies $\mathds{1}_{\mathcal{Q}}\vf \in \mathcal{M}^{10/7,2}_2$ and, since we have for the second parameter defining this Morrey space that $5/3<2<5/2$, we fulfill the condition over $\vf$ stated for the small gradients criterion.
\item[$(iii)$] For the last part of the theorem, we will need even more regularity for the force, indeed, from the previous lines we see that $\vf \in  L^2_t L^2_x$ will not be enough since $2<5/2$. Thus, if we want to work with classical spaces we may ask $\vf\in L^2_tH^1_x$. Indeed, since $ L^2_tH^1_x \subset \mathcal{M}^{2,10/3}_2$ and since $\mathcal{Q}$ is a bounded subset, we find that $\mathds{1}_{\mathcal{Q}}\vf\in L^2_t H^1_x $ implies $\mathds{1}_{\mathcal{Q}}\vf \in\mathcal{M}^{10/7,10/3}_2$. Note in particular that we have here $5/2<10/3$ for the second parameter of the previous Morrey space and thus assuming $\vf \in L^2_tH^1_x$ we satisfy the required hypothesis.
\end{itemize}
%%%%%%%%%%%%%%%%%%%%%%%%%%%%%%%%%%%%%%%%%%%%%%
\section{Proof of the Theorem \ref{Theorem_1}}\label{Secc_Preuve_Theo1}
In this section we will prove Theorem \ref{Theorem_1} with the help of the Theorem \ref{theo_CLM} below for which we  follow the global structure of Theorem \ref{theo_kukavica}, \emph{i.e.} we will decompose each step in function of the hypotheses needed for the force.
%%%%%%%%%%%%%%%%%%%%%%%%%%%%%%%%%%%%%%%%%%%%%%
\begin{Theoreme}\label{theo_CLM}  Let $\Omega$ be a bounded domain of $\mathbb{R}^{3}$ of the form $\Omega=Q_{t,x,r}$ given in (\ref{Def_BouleParabolique}) for some $t\in\mathbb{R}$, $x\in\mathbb{R}^3$ and $\rho>0$. Let $\vu$ be a weak solution for the Navier--Stokes equations (\ref{Equation_NS}) such that $\vu \in L^{\infty}_t L^2_x (\Omega)\cap L^2_t \dot{H}^1_x (\Omega)$, $p\in \mathcal{D}'(\Omega)$, and $\vf \in  \mathcal{D}'(\Omega)$.  Then we have the following points:

\begin{itemize}
 \item[{\bf 1)}] {\bf Energy inequality.} If $\vf\in L^{10/7}_{t}L^{10/7}_{x}(\Omega)$, then the quantity
\begin{equation*}\begin{split}
M=-\partial_{t}|\vu|^{2}+\nu \Delta |\vu|^{2}-2\nu|\vn \otimes \vu|^{2}-&\Div(|\vu|^{2}\vu)\\& -2\langle \Div(p\vu)\rangle+ 2 \vu\cdot \vf,\end{split}
\end{equation*}
is well-defined as a distribution. The solution $\vec u$ is called {\rm dissipative} if $M$ is a locally finite non-negative Borel measure on $\Omega$. 
 \item[{\bf 2)}] {\bf The small gradients criterion.}  Assume that:
\begin{itemize}
\item[$\bullet$]  $\vf\in L^2_tL^2_x(\Omega)$,
\item[$\bullet$]   $\vu$ is dissipative.   
\end{itemize}
There exists positive constants $\epsilon^*>0$ and $\tau_1>5$  which depend only on $\nu$     such that, if $(t_0,x_0)\in \Omega$ and 
 $$\limsup_{r\rightarrow 0} \frac{1}{r}\iint_{]t_0-r^2,t_0+r^2[\times B_{x_0,r}} \vert\vec\nabla\otimes \vec u(s,y)\vert^2\, ds\, dy<\epsilon^*,$$ 
 then there exists a neighborhood $\mathcal{Q}=Q_{t_{0}, x_{0}, \bar{r}}$  of $(t_0,x_0)$  such that $\mathds{1}_{\mathcal{Q}} \,\vec u\in \mathcal{M}^{3,\tau_1}_2$.  
 \item[{\bf 3)}] {\bf Regular points. } Assume that   there exists a    neighborhood $\mathcal{Q}=Q_{t_{0}, x_{0}, \bar{r}}$  of $(t_0,x_0)\in \Omega$  such that 
\begin{itemize}
\item[$\bullet$]  $\mathds{1}_{\mathcal{Q}}\, \vec u\in \mathcal{M}^{3,\tau_1}_2$ for some $\tau_1>5$,
\item[$\bullet$] $\mathds{1}_{\mathcal{Q}} \, \vf \in L^2_tH^1_x$.  
\end{itemize}
Then there exists $r'<\bar{r}$  such that $\vu$ is bounded  on $Q_{t_0,x_0,r'}$. In particular,  the point $(t_0,x_0)$ is regular.
\end{itemize}
 \end{Theoreme}
%%%%%%%%%%%%%%%%%%%%%%%%%%%%%%%%%%%%%%%%%%%%%%
It is worth noting here that at each one of these steps, we have changed the hypotheses for the force $\vf$ (from less regular to more regular) in order to ensure the desired conclusion and we will follow this frame in the proof of the theorem. We do not claim here any kind of optimality with respect to these assumptions.

\begin{Remarque} In the assumptions of Theorem \ref{theo_CLM}, we may add the assumption that $\Div(\vf)=0$.
\end{Remarque} 
Indeed, if $\vf$ belongs to $L^{10/7}_t L^{10/7}_x(\Omega)$, $L^2_t L^2_x(\Omega)$ or $L^2_t H^1_x(\Omega)$, with  $\Omega=Q_{t, x, \rho}$, then it can be extended to $]t-\rho^2,t+\rho^2[\times\mathbb{R}^3$ and still belong to $L^{10/7}_t L^{10/7}_x $, $L^2_t L^2_x $ or $L^2_t H^1_x $. Then, using the fact that the Leray projection operator is bounded on $L^{10/7}_x$, $L^2_x$ and $H^1_x$, we find that we may write $\vf$ as $\vf=\vf_0-\vec\nabla q$ where $\Div(\vf_0)=0$ and $\vf_0$ belongs to $L^{10/7}_t L^{10/7}_x(\Omega)$, $L^2_t L^2_x(\Omega)$ or $L^2_t H^1_x(\Omega)$. Thus we change the couple pressure-force $(p,\vf)$ in the Naver--Stokes equations into $(p+q,\vf_0)$. As $q$ is regular ($q$ and $\vec\nabla q$ belong to $L^{10/7}_t L^{10/7}_x(\Omega)$), we see easily that
$$\underset{\varepsilon \to 0}{\lim}\,\underset{\alpha \to 0}{\lim}\, \Div\left[\big(q\ast \varphi_{\alpha, \varepsilon}\big) \times \big(\vu \ast \varphi_{\alpha, \varepsilon}\big)\right]=\Div(q\vu)=\vu\cdot\vec\nabla q.
$$
In particular, we obtain the same distribution $
M$ when we compute it as associated to the solution $(\vu,p)$ and  the force $\vf$ or as associated to the solution $(\vu,p+q)$ and the divergence-free force $\vf_0$.

\begin{Remarque}\label{Remarque1}Suitability implies dissipativity.
\end{Remarque}
 If we assume a little regularity on $p$ (such as $p\in L^{q_0}_t L^{q_0}_x(\Omega)$ with $q_0>1$), then we shall have $\langle\Div(p\vu)\rangle=  \Div(p\vu)$ and $M=\mu$. In that case, Theorem \ref{theo_CLM} is reduced to Theorem \ref{theo_kukavica}. However, in our theorem, we assume no regularity at all on $p$, so that in particular, $\vec u$ can not be regular in the time variable, as shown by Serrin's counterexample.\\

In fact, the class of dissipative solutions  is \emph{strictly} larger than the class of suitable solutions. It is indeed easy to check that Serrin's counterexample (which is stated without a force) is actually a dissipative solution. 
%%%%%%%
Indeed, we recall that $\vu(t,x)=\phi(t)\vn \psi(x)$ on $]0,1[\times B(0,1)$ where $\phi$ is a bounded function on $\mathbb{R}$ and $\psi$ is a harmonic on $\mathbb{R}^{3}$ and the pressure is given by $p(t,x)=-\frac{|\vu(t,x)|^{2}}{2}-\partial_{t}\phi(t) \psi(x)$. Now we need to verify that the distribution $M$ given in (\ref{Dissipative1}) is non-negative:
$$M=-\partial_{t}|\vu|^{2}+\nu \Delta |\vu|^{2}-2\nu|\vn \otimes \vu|^{2}-\Div(|\vu|^{2}\vu)-2\langle\Div(p\vu)\rangle,
$$
but since $\Delta \vu=0$ we have $\nu \Delta |\vu|^{2}-2\nu|\vn \otimes \vu|^{2}=2\nu \vu \cdot \Delta \vu=0$ and thus we obtain
\begin{eqnarray*}
M&=&-|\vn \psi |^{2}\partial_{t}(\phi(t)^{2})-\Div(|\vu|^{2}\vu)-2\underset{\varepsilon \to 0}{\lim}\underset{\alpha \to 0}{\lim} \Div\big( (p\ast \varphi_{\alpha, \varepsilon}) (\vu \ast \varphi_{\alpha, \varepsilon})\big)\\
&=&-|\vn \psi |^{2}\partial_{t}(\phi(t)^{2})-\Div(|\vu|^{2}\vu)-2\underset{\varepsilon \to 0}{\lim}\underset{\alpha \to 0}{\lim} \Div\left( \left(\bigg[-\frac{|\vu|^{2}}{2}-\partial_{t}\phi\psi\bigg]\ast \varphi_{\alpha, \varepsilon}\right) (\vu \ast \varphi_{\alpha, \varepsilon})\right)\\
&=&-|\vn \psi |^{2}\partial_{t}(\phi(t)^{2})-\Div(|\vu|^{2}\vu)+2\underset{\varepsilon \to 0}{\lim}\underset{\alpha \to 0}{\lim} \Div\left( \left(\frac{|\vu|^{2}}{2}\ast \varphi_{\alpha, \varepsilon}\right) (\vu \ast \varphi_{\alpha, \varepsilon})\right)\\
&&+ 2\underset{\varepsilon \to 0}{\lim}\underset{\alpha \to 0}{\lim} \Div\left( \left([\partial_{t}\phi \psi]\ast \varphi_{\alpha, \varepsilon}\right) (\phi(t)\vn \psi(x) \ast \varphi_{\alpha, \varepsilon})\right).
\end{eqnarray*}
Recalling that $\Delta \psi=0$ and passing to the limit $\alpha, \varepsilon\longrightarrow 0$ it is easy to see that $M=0$; thus the example of Serrin is dissipative in the sense of the Definition \ref{Definition_Dissipative} and is in the scope of Theorem \ref{Theorem_1}. 
%%%%%%%%%%%%%%%%%%%%%%%%%%%%%%%%%%%%%%%%%%%%%%
\subsection{The new variable}

We start the proof of Theorem \ref{theo_CLM} by considering a function  $\vu$ and a distribution $p$ that satisfy the Navier--Stokes equations (\ref{Equation_NS}) over a bounded domain $\Omega\subset \mathbb{R}\times \mathbb{R}^{3}$, for a divergence-free force $\vf$. In order to simplify the notation, and with no loss of generality (as we are interested in local properties), we will assume once and for all that the set $\Omega$ is of the form 
\begin{equation}\label{Set_Definition}
\Omega=I\times B_{x_{0},\rho},
\end{equation}
where $I=]a,b[$ is an interval and $B_{x_{0},\rho}=B(x_{0}, \rho)$ is an open ball in $\mathbb{R}^{3}$ of radius $\rho >0$ and center $x_{0}\in \mathbb{R}^{3}$. In this section, we only assume that $\vu\in L^\infty_t L^2_x\cap L^2_t \dot{H}^1_x$ and that $\Div(\vf)=0$.\\

Our first step is to consider the curl of $\vu$, which will be denoted by $\vor=\vn \wedge \vu$. We  obtain the following equation:
\begin{eqnarray}\label{Equation_Vorticity_1}
\partial_t \vor=\nu \Delta \vor - \vn \wedge (\vor \wedge \vu)+\vn \wedge \vf,
\end{eqnarray}
where the pressure $p$ has disappeared (since we have $\vn \wedge \vn p \equiv 0$).\\

However, as said in the introduction, we shall be interested in the more regular distribution $\vec u$ than in the distribution $\vec\omega$ and for this we proceed as follows: since we want to study the regularity of $\vu$ inside $\Omega$; we shall restrict ourselves to a smaller domain 
\begin{equation}\label{Def_Omega0}
\Omega_0=I_0\times B_{x_{0},\rho_0},
\end{equation}
with $I_0=]a_0,b_0[$, where $a<a_0<b_0<b$ and $0<\rho_0<\rho$. Then, to go back to $\vu$ from the vorticity $\vec \omega$, we introduce a cut-off function $\psi\in\mathcal{D}(\mathbb{R}\times\mathbb{R}^3)$ which is equal to $1$ on a neighborhood of $\Omega_0$ and is compactly supported in $\Omega$.  More precisely, we ask $\psi$ to be of the form 
$$\psi(t,x)=\phi(t)\Phi(x),$$ 
where $\phi$ is equal to $1$ on a neighborhood of $I_0$ and is compactly supported within $I$, while $\Phi$ is equal to $1$ on a neighborhood of $B_{x_{0},\rho_0}$ and is compactly supported within $B_{x_{0},\rho}$. The distribution $\psi \vec\omega$ may be viewed as defined on the whole $\mathbb{R}\times\mathbb{R}^3$ and clearly belongs to $L^\infty_{t} H_{x}^{-1}\cap L_{t}^2 L_{x}^2$.
Thus,  using this localization function $\psi$, we can define a new function $\vv$ in the following way
\begin{equation}\label{Definition_v}
\vv = -\frac{1}{\Delta} \vn \wedge (\psi \vor)=-\frac{1}{\Delta} \vn \wedge (\psi \vn \wedge \vu).
\end{equation}
Note in particular that, on $\Omega_0$, the derivatives of $\psi$ are equal to $0$, so that
$$ \Delta\vv =  -\vn \wedge (\psi \vor)= - \vn \wedge (  \vn \wedge \vu)=\Delta\vu,$$
(since $\Div(\vu)=0$). We can see then from this identity that, on $\Omega_0$, $\vec v$ is equal to $\vu$ up to a harmonic (in the space variable) correction $\vec w$. Throughout the paper, our stategy will be to replace the study of the regularity of $\vu$ with the study of the regularity of $\vv$, and to link those two regularities by a precise study of the harmonic correction $\vec w=\vec v-\vu$.\\

We begin with some elementary facts on $\vv$ and $\vec w$ :

%%%%%%%%%%%%%%%%%%%%%%%%%%%%%%%%%%%%%%%%%%%%%%
\begin{Proposition}\label{Proposition_1} Under the assumptions   $\vu \in L^\infty_t L^2_x(\Omega)\cap L^2_t\dot H^1_x(\Omega)$ on the solution $\vu$ of  the Navier--Stokes equations (\ref{Equation_NS}),
the function $\vv$ defined by the formula (\ref{Definition_v}) above satisfies the following points \begin{itemize}
\item[1)]    $\Div(\vv)=0$,
\item[2)] $\vv\in L^\infty_t L^2_x(\Omega_0)\cap L^2_t H^1_x(\Omega_0)$,
\item[3)]  the function $\vec w=\vv-\vu$ satisfies $\vec{w}\in L^\infty_t {\rm Lip}_x(\Omega_0)$.
\end{itemize}
\end{Proposition}
%%%%%%%%%%%%%%%%%%%%%%%%%%%%%%%%%%%%%%%%%%%%%%
{\bf \textit{Proof.}} The first point is obvious, since the divergence of a curl is always null. For the second point, we will use the identity 
$$\psi \vn \wedge \vu=\vn \wedge(\psi \vu)-(\vn \psi) \wedge \vu,$$
and using the definition of $\vv$ given above in (\ref{Definition_v}) we obtain the expression
\begin{equation}\label{Formule_utile}
\vv = -\frac{1}{\Delta} \vn \wedge \left(\vn \wedge(\psi \vu)-(\vn \psi) \wedge \vu\right).
\end{equation}
We will prove in the following items that $\vv \in L^\infty_{t}L^2_{x}(\Omega_{0})$ and  $\vv\in L^2_{t}H_{x}^{1}(\Omega_{0})$.

\begin{itemize}
%%%%%%%%%%%%%
\item[$\bullet$] We start with $\vv \in L^\infty_{t}L^2_{x}(\Omega_{0})$. Taking the norm $L^2(B_{x_{0},\rho_0})$ in the space variable of the expression (\ref{Formule_utile}) we have
\begin{eqnarray}
\|\vv\|_{L^{2}(B_{x_{0},\rho_0})}&\leq &\left\| \frac{1}{\Delta} \vn \wedge \left(\vn \wedge(\psi \vu)\right)\right\|_{L^{2}(B_{x_{0},\rho_0})}+ \left\|\frac{1}{\Delta} \vn \wedge \left((\vn \psi) \wedge \vu\right)\right\|_{L^{2}(B_{x_{0},\rho_0})}\nonumber\\
&\leq &\left\|\frac{1}{\Delta} \vn \wedge \big( \vn \wedge(\psi \vu) \big)\right\|_{L^2(\mathbb{R}^{3})}+ \left\|\frac{1}{\Delta} \vn \wedge \big((\vn \psi) \wedge \vu\big)\right\|_{L^2(\mathbb{R}^{3})}.\label{NormeL2_Utile}
\end{eqnarray}
Applying Hardy-Littlewood-Sobolev inequalities in the second term above we obtain
$$\|\vv\|_{L^{2}(B_{x_{0},\rho_0})} \leq C\left\|\psi \vu\right\|_{L^2(\mathbb{R}^{3})}+ C\left\|\vn \psi \wedge \vu\right\|_{L^{6/5}(\mathbb{R}^{3})}.$$
Now, using the support properties of the function $\psi$ and using the H\"older inequality we can write
\begin{equation}\label{NormeL2_Utile1}
\|\vv\|_{L^{2}(B_{x_{0},\rho_0})}  \leq C\|\psi\|_{L^\infty(B_{x_{0},\rho})}\|\vu\|_{L^2(B_{x_{0},\rho})}+C\|\vn \psi\|_{L^{3}(B_{x_{0},\rho})}\|\vu\|_{L^2(B_{x_{0},\rho})}.
\end{equation}
It remains to take the $L^{\infty}$ norm in the time variable in order to obtain
\begin{eqnarray*}
\|\vv\|_{L^{\infty}(I_0,L^{2}(B_{x_{0},\rho_0}))} &\leq &C\left(\|\psi\|_{L^{\infty}_tL^\infty_x}+\|\vn \psi\|_{L^{\infty}_t L^3_x}\right)\|\vu\|_{L^{\infty}(I_0,L^2(B_{x_{0},\rho}))}\\
&\leq & C_{\psi} \|\vu\|_{L^{\infty}(I_0,L^2(B_{x_{0},\rho}))}\leq  C_{\psi} \|\vu\|_{L_{t}^{\infty}L_{x}^2(\Omega)},
\end{eqnarray*}
since $ I_0\times B_{x_{0},\rho}  \subset \Omega$. The last quantity above is bounded by the hypotheses on $\vu$. 
%%%%%%%%%%%%%
\item[$\bullet$] We study now the fact that $\vv\in L^2\big(I_0, H^{1}(B_{x_{0},\rho_0})\big)$: taking the $H^{1}(B_{x_{0},\rho_0})$ norm in the expression (\ref{Formule_utile}) we obtain
\begin{eqnarray*}
\|\vv\|_{H^{1}(B_{x_{0},\rho_0})} &\leq &\left\|\frac{1}{\Delta} \vn \wedge \big( \vn \wedge(\psi \vu) \big)\right\|_{H^1(\mathbb{R}^{3})}+ \left\|\frac{1}{\Delta} \vn \wedge \big((\vn \psi) \wedge \vu\big)\right\|_{H^1(\mathbb{R}^{3})}.
\end{eqnarray*}
The first term above can be controlled by $\|\psi \vu\|_{H^1(\mathbb{R}^{3})}$, thus (due to the support properties of the function $\psi$) by $C_{\psi}  \|\vu\|_{H^1(B_{x_{0},\rho})}$. For the second term of the previous expression, we have by definition of the $H^{1}$ norm
\begin{equation*}\begin{split} \left\|\frac{1}{\Delta} \vn \wedge \big((\vn \psi) \wedge \vu\big)\right\|_{H^1(\mathbb{R}^{3})}=&  \left\|\frac{1}{\Delta} \vn \wedge \big((\vn \psi) \wedge \vu\big)\right\|_{L^{2}(\mathbb{R}^{3})}\\ &+\left\|\vn\left(\frac{1}{\Delta} \vn \wedge \big((\vn \psi) \wedge \vu\big)\right)\right\|_{L^{2}(\mathbb{R}^{3})}.
\end{split}\end{equation*}
Following the same computations performed in (\ref{NormeL2_Utile})-(\ref{NormeL2_Utile1}) we see that the first quantity in the right-hand side of the previous formula is controlled by $\|\vn \psi\|_{L^{3}(B_{x_{0},\rho})} \|\vu\|_{L^2(B_{x_{0},\rho})}$, while the second quantity in the right-hand side can be estimated by  $\|\vn\psi \|_{L^{\infty}(B_{x_{0},\rho})} \|\vu\|_{L^2(B_{x_{0},\rho})}$. Gathering all this estimates we obtain
$$\|\vv\|_{H^{1}(B_{x_{0},\rho_0})}\leq C_{\psi}\|\vu\|_{H^1(B_{x_{0},\rho})},$$
and thus we have
$$\|\vv\|_{L^{2}(I_0,H^{1}(B_{x_{0},\rho_0}))}\leq C_{\psi}\|\vu\|_{L^{2}(I_0,H^1(B_{x_{0},\rho}))}\leq  C_{\psi}\|\vu\|_{L_{t}^{2} H_{x}^1(\Omega)}<+\infty.$$
%%%%%%%%%%%%%
\end{itemize}$\ $\\

We now prove the third point of the proposition. For this, we start using the following general identities (where we use $\Div(\psi\vu)=(\vu\cdot\vn\psi)$, as $\Div(\vu)=0$)
\begin{eqnarray*}
\vv &=& -\frac{1}{\Delta}\vn \wedge \left(\vn \wedge (\psi \vu)-(\vn \psi) \wedge \vu\right)
\\ &=& -\frac{1}{\Delta}\left[\vn \wedge \left(\vn \wedge (\psi \vu)\right)\right] + \frac{1}{\Delta}\left[ \vn \wedge \left((\vn \psi) \wedge \vu\right)\right]\\
&=&-\frac{1}{\Delta}\left[\vn (\vu \cdot \vn \psi) - \Delta (\psi\vu)\right]+ \frac{1}{\Delta}\left[ \vn \wedge \left((\vn \psi) \wedge \vu\right)\right] \\\ &=&-\frac{1}{\Delta}\left[\vn (\vu \cdot \vn \psi)\right]+ \psi\vu + \frac{1}{\Delta}\left[ \vn \wedge \left((\vn \psi) \wedge \vu\right)\right].
\end{eqnarray*}
From this last identity, it is possible to derive a reformulation for $\vv$:
\begin{equation}\label{ReDefinition_Vu}
\vv=\psi\vu+ \frac{1}{\Delta}\left[ \vn \wedge \left((\vn \psi) \wedge \vu\right)\right]-\frac{1}{\Delta}\left[\vn (\vu \cdot \vn \psi)\right].
\end{equation}
Now, since by definition we have that $\psi\equiv 1$ over $\Omega_0$, we obtain  the following decomposition on $\Omega_0$
$$\vv=\vu+ \vec{w},$$
where 
\begin{equation}\label{Def_DifferenceUV}
\vec{w}=\frac{1}{\Delta}\left[ \vn \wedge \left((\vn \psi) \wedge \vu\right)\right]-\frac{1}{\Delta}\left[\vn (\vu \cdot \vn \psi)\right].  
\end{equation}
We recall now that the operator $\frac{1}{\Delta}$ is given by convolution with a kernel $K$: indeed, for an admissible function $f$ we have
\begin{equation}\label{Def_Kernel}
\frac{1}{\Delta}f(x)=K\ast f(x)=-\frac{1}{4\pi}\int_{\mathbb{R}^{3}}\frac{f(y)}{|x-y|}dy.
\end{equation}
For $(t,x)\in \Omega_0$, we have $x\in B_{x_{0},\rho_0}$ while $\psi(t,y)=1$ on $I_0\times B_{x_{0},\rho_1}$ for some $\rho_1$ with $\rho_0<\rho_1<\rho$. In particular, $\vec \nabla\psi(t,y)$ is  identically null for $t\in I_0$, $x\in B_{x_{0},\rho_0}$ and $\vert y-x\vert<\rho_1-\rho_0$. Thus, with the definition of $\vw$ given in (\ref{Def_DifferenceUV}) above we may write, for every multi-index 
 $\alpha$ in $\mathbb{N}^{3}$  
\begin{eqnarray*}
|\partial_{x}^{\alpha}\vec{ w}(t,x)|&\leq&\left|(\partial_{x}^{\alpha}K)\ast \left[ \vn \wedge \left((\vn \psi) \wedge \vu\right)\right](t,x)\right|+\left|(\partial_{x}^{\alpha}K)\ast\left[\vn (\vu \cdot \vn \psi)\right](t,x)\right|\\
&\leq&\sum_{i,j,k}^{3}\left|(\partial_{x}^{\alpha}K)\ast \left[ \partial_{x_{i}} \left((\partial_{x_{j}}\psi) u_{k}\right)\right](t,x)\right|+\sum_{i,j,k}^{3}\left|(\partial_{x}^{\alpha}K)\ast\left[\partial_{x_{i}} (u_{j} \partial_{x_{k}}\psi)\right](t,x)\right|\\
&\leq&\sum_{i,j,k}^{3}\left|(\partial_{x}^{\alpha}\partial_{x_{i}}K)\ast \left[(\partial_{x_{j}}\psi) u_{k}\right](t,x)\right|+\sum_{i,j,k}^{3}\left|(\partial_{x}^{\alpha}\partial_{x_{i}}K)\ast\left[ u_{j} (\partial_{x_{k}}\psi)\right](t,x)\right|,
\end{eqnarray*}
and obtain
\begin{equation}\label{Formule_Kernel}
|\partial_{x}^{\alpha}\vec{ w}(t,x)|\leq C_{\alpha}\int_{\{\rho_1-\rho_0<\vert y-x\vert, \ y\in B_{x_{0},\rho} \}}\frac{|\vn \psi(t,y)| |\vu (t,y)|}{|x-y|^{2+|\alpha|}}dy, 
\end{equation}
Thus, we have the following control for $(t,x)\in\Omega_0$ :
$$|\partial_{x}^{\alpha}\vec{ w}(t,x)|\leq \frac{C}{(\rho_{1}-\rho_{0})^{2+|\alpha|}} \|\vn \psi(t,\cdot)\|_{L^{2}(B_{x_{0},\rho})} \|\vu(t,\cdot)\|_{L^{2}(B_{x_{0},\rho})},$$
from which we obtain that $\partial_{x}^{\alpha}\vec{w}\in L^{\infty}_{t}L^{\infty}_{x}(\Omega_0)$. The Proposition \ref{Proposition_1} is now completely proven.  \hfill $\blacksquare$\\
\begin{Remarque}\label{Remarque3} In Proposition \ref{Proposition_1} we have stated the results over the set $\Omega_{0}\subset \Omega$ and it is possible to extend some properties of the new variable $\vv$ to the set $\Omega$. Indeed, following the same computations performed in the second point of the previous proposition it is easy to see that we have 
$$\vv\in L^\infty_t L^2_x(\Omega)\cap L^2_t H^1_x(\Omega).$$  
\end{Remarque}
However, the fact $\vw\in L^\infty_t {\rm Lip}_x$ can not be extended to the set $\Omega$ and in order to obtain this property we need to work over a smaller subset $\Omega_{0}\subset \Omega$. Remark nevertheless that over $\Omega$ we have:
%%%%%%%%%%%%%%%%%%%%%%%%%%%%%%%%%%%%%%%%%%%%%%
\begin{Corollaire}\label{Proprietes_DifferenceUV}
If we work over all the subset $\Omega$ since $\vu\in L^\infty_t L^2_x(\Omega)\cap L^2_t H^1_x(\Omega)$ and $\vv\in L^\infty_t L^2_x(\Omega)\cap L^2_t H^1_x(\Omega)$, we also obtain $\vw\in L^\infty_t L^2_x(\Omega)\cap L^2_t H^1_x(\Omega)$. Indeed we have
$$\|\vw\|_{L_{t}^{\infty} L_{x}^2(\Omega)}\leq C_{\rho, \psi}\|\vu\|_{L_{t}^{\infty} L_{x}^2(\Omega)} \quad \mbox{and}\quad \|\vw\|_{L_{t}^{2} H_{x}^1(\Omega)}\leq C_{\rho, \psi}\|\vu\|_{L_{t}^{2} H_{x}^1(\Omega)}.$$
\end{Corollaire}
%%%%%%%%%%%%%%%%%%%%%%%%%%%%%%%%%%%%%%%%%%%%%%
\subsection{The companion equation.}
Now we turn our attention to the equation satisfied by the new variable $\vv$ and we will see here that we obtain a new Navier--Stokes equation on $\Omega_0$
\begin{equation*}
\partial_t\vv=\nu\Delta\vv-(\vv\cdot\vn)\vv-\vn q+\vec F,
\end{equation*}
where $\vn q$ is a gradient term and $\vec F$ is a divergence-free force. As said in the introduction, this equation will be called the \emph{companion equation} of the original Navier--Stokes equations (\ref{Equation_NS}). Of course, our aim is now to prove that $q$ and $\vec F$ may be easily estimated from the assumptions on $\vf$ and $\vu$ without involving the pressure $p$.\\ 

We recall that the variable $\vv$ was defined by $\vv = -\frac{1}{\Delta} \vn \wedge (\psi \vn \wedge \vu)$, we shall thus define now in an analogous way
 \begin{equation*} 
\vec F_0 = -\frac{1}{\Delta} \vn \wedge (\psi \vn \wedge \vf),
\end{equation*}
and we have the following lemma which is a consequence of the previous computations performed in Proposition \ref{Proposition_1}.
%%%%%%%%%%%%%%%%%%%%%%%%%%%%%%%%%%%%%%%%%%%%%%
\begin{Lemme}\label{LemmeForce0}
Let $\Omega$ be a bounded subset of $\mathbb{R}\times \mathbb{R}^{3}$ of the form (\ref{Set_Definition}) and let $\Omega_{0}$ be the set given in (\ref{Def_Omega0}). Let $\vf$ be a given force such that $\Div(\vf)=0$, then
\begin{itemize}
\item[1)] if $\vf \in L^{10/7}_{t}L^{10/7}_{x}(\Omega)$ then $\vF_{0} \in L^{10/7}_{t}L^{10/7}_{x}(\Omega_{0})$,
\item[2)] if $\vf \in L^{2}_{t}L^{2}_{x}(\Omega)$ then $\vF_{0} \in L^{2}_{t}L^{2}_{x}(\Omega_{0})$,
\item[3)] if $\vf \in L^{2}_{t}H^{1}_{x}(\Omega)$ then $\vF_{0} \in L^{2}_{t}H^{1}_{x}(\Omega_{0})$.
\end{itemize}
\end{Lemme}
%%%%%%%%%%%%%%%%%%%%%%%%%%%%%%%%%%%%%%%%%%%%%%
{\bf \textit{Proof.}} It is enough to remark that, in the same spirit of formulas  (\ref{ReDefinition_Vu})-(\ref{Def_DifferenceUV}), the force $\vF_{0}$ can be decomposed over $\Omega_{0}$ by $\vF_{0}=\vf+\vw_{f}$ where 
$$\vw_{f}= \frac{1}{\Delta}\left[ \vn \wedge \left((\vn \psi) \wedge \vf\right)\right]-\frac{1}{\Delta}\left[\vn (\vf \cdot \vn \psi)\right],$$
and $\vw_{f}\in L^\infty_t {\rm Lip}_x(\Omega_{0})$. \hfill $\blacksquare$\\

Remark also that, as the distribution $ \vn \wedge (\psi \vn \wedge \vf)$ is compactly supported and as $ -\frac{1}{\Delta} $ is a convolution operator with the distribution $\delta_t\otimes \frac{1}{4\pi| x|}$, the quantity $\vec F_0$ is well defined for any distribution $\vf\in\mathcal{D}'(\Omega)$: we can consider a wider framework for the force and this point of view will be adopted in the following proposition.\\ 

Indeed, we will explain with the next result how the deduce the companion equation satisfied by the variable $\vv$ and we will study the relationship between $\vec{F}_{0}$ and the force $\vec{F}$. We will also give some important properties of the new pressure $q$. 
%%%%%%%%%%%%%%%%%%%%%%%%%%%%%%%%%%%%%%%%%%%%%%
\begin{Proposition}\label{Proposition_2} Under the assumptions   $\vu \in L^\infty_t L^2_x\cap L^2_t\dot H^1_x(\Omega)$, $p\in\mathcal{D}'(\Omega)$ and $\vf\in\mathcal{D}'(\Omega)$,  on the solution $\vu$ of  the Navier--Stokes equations (\ref{Equation_NS}),
the function $\vv$ defined by the formula (\ref{Definition_v}) above
 satisfies the following Navier--Stokes equations on $\Omega_0$
\begin{equation}\label{Equation_V0}
\partial_t \vv =\nu \Delta \vv -  (\vv \cdot\vn) \vv- \vn q+\vF,
\end{equation}
with the following properties
\begin{itemize}
\item[1)]  
the pressure $q$ is a function such that $q\in L^{3/2}_tL^{3/2}_x(\Omega_0)$,
\item[2)] the force $\vF$ is such that $\Div(\vF)=0$ and $\vF-\vec F_0\in L^2_tL^2_x( \Omega_0)$, where $\vec F_0 = -\frac{1}{\Delta} \vn \wedge (\psi \vn \wedge \vf)$.
\end{itemize}
\end{Proposition}
Observe at this stage of the proof that we only assume that $p\in\mathcal{D}'(\Omega)$ and $\vf\in\mathcal{D}'(\Omega)$, but this is enough to obtain that the new pressure $q$ belongs to the space $L^{3/2}_tL^{3/2}_x(\Omega_0)$. However we will need later on some extra assumptions on $\vf$ in order to obtain a more regular behavior for the global force $\vF$.\\

%%%%%%%%%%%%%%%%%%%%%%%%%%%%%%%%%%%%%%%%%%%%%%
{\bf \textit{Proof.}} We start by describing the equation satisfied by $\partial_{t}\vv$. Since we are working on $\Omega_0$, we have $\partial_{t}\psi=0$ for $t\in I_0$ by the support properties of $\psi$ and we may write:
$$\partial_t \vv=\partial_t \left[-\frac{1}{\Delta} \vn \wedge (\psi \vor)\right]=-\frac{1}{\Delta} \vn \wedge \big(\psi (\partial_t \vor)\big),$$
thus, using equation (\ref{Equation_Vorticity_1}) we obtain
\begin{eqnarray}
\partial_t \vv &=& -\frac{1}{\Delta} \vn \wedge \left(\psi\left(\nu\Delta \vor -  \vn \wedge (\vor \wedge \vu)+\vn \wedge \vf\right)\right)\nonumber\\
&=& \nu\left(\underbrace{-\frac{1}{\Delta} \vn \wedge \left(\psi\Delta \vor \right)}_{(A)}\right) + \underbrace{\frac{1}{\Delta} \vn \wedge \left( \psi (\vn \wedge (\vor \wedge \vu))\right)}_{(B)}+ \vF_0,\label{Equation_V_1}
\end{eqnarray}
where 
\begin{equation}\label{Definition_F0} 
\vF_0=- \frac{1}{\Delta} \vn \wedge \left(\psi \vn \wedge \vf\right).
\end{equation}
We study now each one of the terms $(A)$ and $(B)$ in order to simplify the expression (\ref{Equation_V_1}). 
\begin{itemize}
\item[(A)] For the first term $-\frac{1}{\Delta} \vn \wedge \left(\psi\Delta \vor \right)$ we write the following identities for $\psi\Delta \vor$:
\begin{eqnarray*}
 \psi \Delta \vor = \psi \Delta (\vn \wedge \vu) &= &\Delta (\psi \vn \wedge \vu) - (\Delta \psi) \vn \wedge \vu -2 \sum_{j=1}^3 (\partial_{x_j} \psi)(\partial_{x_j} \vn \wedge \vu)\\
&=& \Delta (\psi \vn \wedge \vu) - (\Delta \psi) \vn \wedge \vu -2 \sum_{j=1}^3\partial_{x_j}\left( (\partial_{x_j} \psi)( \vn \wedge \vu)\right)\\&& +2 (\Delta \psi) \vn \wedge \vu\\
&=& \Delta (\psi \vn \wedge \vu) + (\Delta \psi) \vn \wedge \vu -2 \sum_{j=1}^3\partial_{x_j}\left( (\partial_{x_j} \psi)( \vn \wedge \vu)\right).
\end{eqnarray*}
Now, using the classical vector calculus identitiy 
$$ \vn\wedge(a\vec b)=a\vn\wedge\vec b+(\vn a)\wedge\vec b,$$ we obtain
\begin{eqnarray*}
 \psi \Delta \vor&=& \Delta (\psi \vn \wedge \vu) + \left(\vn \wedge(\Delta\psi \vu) - \vn (\Delta \psi)\wedge \vu\right) -2 \sum_{j=1}^3\partial_{x_j}\left(  \vn \wedge\left( (\partial_{x_j} \psi)  \vu\right)\right)\\ &&+2 \sum_{j=1}^3\partial_{x_j}\left(\vn (\partial_{x_j}\psi) \wedge \vu\right),
\end{eqnarray*}
so that we have
\begin{eqnarray*}
-\frac{1}{\Delta} \vn \wedge[\psi \Delta \vor]&= &\underbrace{-\frac{1}{\Delta} \vn \wedge\left[\Delta (\psi \vn \wedge \vu)\right]}_{(1)} -\underbrace{\frac{1}{\Delta} \vn \wedge\left[\vn \wedge(\Delta\psi \vu)- \vn (\Delta \psi)\wedge \vu\right]}_{(2)} \\
& & + \underbrace{2 \frac{1}{\Delta} \vn \wedge\left[ \sum_{j=1}^3\partial_{x_j}\left(  \vn \wedge\left( (\partial_{x_j} \psi)  \vu\right)\right)\right]}_{(3)} \\ & &- \underbrace{2\frac{1}{\Delta} \vn \wedge \left[\sum_{j=1}^3\partial_{x_j}\left(\vn (\partial_{x_j}\psi) \wedge \vu\right)\right]}_{(4)}.
\end{eqnarray*}
At this point we remark that the term $(1)$ above is in fact $\Delta \vv$, indeed:
$$-\frac{1}{\Delta} \vn \wedge\left[\Delta (\psi \vn \wedge \vu)\right] =\Delta  \left(-\frac{1}{\Delta} \vn \wedge(\psi \vn \wedge \vu)\right)=\Delta \vv,$$
and then we can write 
\begin{equation*}
-\frac{1}{\Delta} \vn \wedge(\psi \Delta \vor)=\Delta \vv + \vF_1,
\end{equation*}
where $\vF_{1}=(2)+(3)+(4)$, \emph{i.e.}:
\begin{eqnarray}
\vF_1 &= &-\frac{1}{\Delta} \vn \wedge\left[\vn \wedge(\Delta\psi \vu) - \vn (\Delta \psi)\wedge \vu\right] + 2 \frac{1}{\Delta} \vn \wedge\left[ \sum_{j=1}^3\partial_{x_j}\left(  \vn \wedge\left( (\partial_{x_j} \psi)  \vu\right)\right)\right]\nonumber \\
& &- 2\frac{1}{\Delta} \vn \wedge \left[\sum_{j=1}^3\partial_{x_j}\left(\vn (\partial_{x_j}\psi) \wedge \vu\right)\right].\label{Definition_F1}
\end{eqnarray}
Observe that since the divergence of a curl is always null, we obtain that $\Div(\vF_{1})=0$.
\item[(B)] We study now the second term of (\ref{Equation_V_1}). For $ \psi (\vn \wedge (\vor \wedge \vu))$, using vector calculus identities we write
$$\psi (\vn \wedge (\vor \wedge \vu))=\vn \wedge  \psi(\vor \wedge\vu ) -(\vn \psi)\wedge (\vor \wedge \vu)=\vn \wedge (\vor \wedge \psi\vu ) -(\vn \psi)\wedge (\vor \wedge \vu),$$
so we obtain
$$\frac{1}{\Delta} \vn \wedge \left[\psi (\vn \wedge (\vor \wedge \vu))\right]=\frac{1}{\Delta} \vn \wedge \left[\vn \wedge (\vor \wedge \psi\vu )\right]-\frac{1}{\Delta} \vn \wedge \left[(\vn \psi)\wedge (\vor \wedge \vu)\right].$$
We remark here that 
$$\vn \wedge \left[\vn \wedge (\vor \wedge \psi\vu )\right]=\left(\vn (\Div (\vor \wedge \psi\vu )) - \Delta (\vor \wedge \psi\vu) \right)$$
and we can write
\begin{eqnarray*}
\frac{1}{\Delta} \vn \wedge \left[\psi (\vn \wedge (\vor \wedge \vu))\right]& =& \frac{1}{\Delta}\left(\vn (\Div (\vor \wedge \psi\vu )) - \Delta (\vor \wedge \psi\vu) \right) -\frac{1}{\Delta} \vn \wedge \left[(\vn \psi)\wedge (\vor \wedge \vu)\right]\\
&=&  -\vn \left(-\frac{1}{\Delta}(\Div (\vor \wedge \psi\vu ) \right)  -\vor \wedge \psi\vu  -\frac{1}{\Delta} \vn \wedge \left[(\vn \psi)\wedge (\vor \wedge \vu)\right].
\end{eqnarray*}
Finally, the second term of  (\ref{Equation_V_1})  can be rewritten in the following form
\begin{equation*}
\frac{1}{\Delta} \vn \wedge \left[\psi (\vn \wedge (\vor \wedge \vu))\right]= -\vn q_1 -\vor \wedge \psi\vu + \vF_2.
\end{equation*}
where we have
\begin{equation}\label{Definition_q1}
 q_1=-\frac{1}{\Delta}(\vn \cdot (\vor \wedge \psi\vu )),
\end{equation}
and 
\begin{equation}\label{Definition_F2}
\vF_2=-\frac{1}{\Delta} \vn \wedge \left[(\vn \psi)\wedge (\vor \wedge \vu)\right].
\end{equation}
\end{itemize}
Remark in particular that we have $\Div(\vec F_2)=0$.\\

Now, coming back to the equation (\ref{Equation_V_1}) and with the definition of the quantities $q_{1}, \vF_{0}, \vF_{1}$ and $\vF_{2}$ given in  (\ref{Definition_q1}), (\ref{Definition_F0}), (\ref{Definition_F1}) and  (\ref{Definition_F2}) respectively, we obtain the following equation
\begin{equation}\label{Equation_Intermediaire_1}
\partial_t \vv=\nu\Delta \vv -\big[\vor \wedge \psi\vu\big] -\vn q_1+ \vF_0 + \nu\vF_1 + \vF_2 ,
\end{equation}
which is almost the desired equation (\ref{Equation_V0}) stated in Proposition \ref{Proposition_2}, but we still need to study the term $\vor \wedge \psi\vu$. For this,  recalling that on $\Omega_0$, the function $\psi$ is \emph{constant} and equal to one,   we have the identity $\vor=\vn \wedge \psi \vu$ and we can write
\begin{equation*}
\vor \wedge \psi \vu= \big(\vn \wedge \psi \vu\big)\; \wedge \psi \vu.\end{equation*}
Now, we rewrite  $\psi \vu$ using the formula (\ref{ReDefinition_Vu}) given in page \pageref{ReDefinition_Vu}:
$$\psi \vu =\vv +\frac{1}{\Delta}\left[\vn (\vu \cdot \vn \psi)\right]- \frac{1}{\Delta}\left[\vn \wedge \left((\vn \psi) \wedge \vu\right)\right] =\vv+\vec\beta,
$$
where we have defined $\vec\beta$ by the expression
\begin{equation}\label{Definition_Beta}
\vec\beta=\frac{1}{\Delta}\left[\vn (\vu \cdot \vn \psi)\right]- \frac{1}{\Delta}\left[\vn \wedge \left((\vn \psi) \wedge \vu\right)\right].
\end{equation}
\begin{Remarque}
Observe that we have $\vec{\beta}=-\vw$, where $\vw$ was given in (\ref{Def_DifferenceUV}). 
\end{Remarque}

We thus obtain the formula
\begin{eqnarray}
\vor \wedge \psi \vu&= &\vn \wedge (\vv+\vec \beta) \wedge (\vv +\vec \beta)\nonumber \\
&=&(\vn \wedge \vv) \wedge \vv + (\vn \wedge \vv) \wedge\vec \beta + (\vn\wedge\vec\beta) \wedge \vv +  (\vn\wedge\vec\beta) \wedge \vec\beta.\label{Identite_Beta}
\end{eqnarray}
Since the new variable $\vv$ is divergence free (by Proposition \ref{Proposition_1}) we have the identity $(\vn \wedge \vv) \wedge \vv= (\vv\cdot\vn)  \vv-\frac{1}{2}\vn |\vv|^{2}$. 
Let us define now $q_{3}$ by
\begin{equation}\label{definition_q3}
q_{3}=-\frac{1}{2}|\vv|^{2},
\end{equation}
and then we have for the first term of (\ref{Identite_Beta}):
$$(\vn \wedge \vv) \wedge \vv=(\vv\cdot\vn)  \vv+\vn q_{3}.$$
In order to estimate the remaining therms of (\ref{Identite_Beta}) we define 
\begin{equation}\label{Definition_A}
\vec A= (\vn \wedge \vv) \wedge\vec \beta +  (\vn\wedge\vec\beta) \wedge \vv  + (\vn\wedge\vec\beta) \wedge \vec\beta,
\end{equation} 
we remark now that on $\Omega_0$, we have $\vec A=\psi \vec A$ and we decompose $\psi\vec{A}$ in the following manner:
\begin{equation}\label{formule_vorticite2}
\psi \vec{A}=-\vF_{3}+\vn q_{2},
\end{equation}
where 
\begin{equation}\label{definition_q2}
q_{2}=\frac{1}{\Delta} \Div(\psi \vec{A}). 
\end{equation}
Again, remark that we have $\Div(\vF_3)=0$.\\

Thus, we have obtained that $\vor \wedge \psi \vu= (\vv\cdot\vn) \vv - \vF_{3}+\vn q_{2}+\vn q_{3}$ and getting back to (\ref{Equation_Intermediaire_1}) we can write
\begin{eqnarray*}
\partial_t \vv&=&\nu\Delta \vv -\big[\vor \wedge \psi\vu\big] -\vn q_1+ \vF_0 + \nu\vF_1 + \vF_2\\
&=& \nu\Delta \vv -\left[  (\vv\cdot\vn) \vv - \vF_{3}+\vn q_{2}+\vn q_{3}\right] -\vn q_1+ \vF_0 + \nu\vF_1 + \vF_2,
\end{eqnarray*} 
and finally we obtain the companion equation for $\vv$
\begin{equation*}
\partial_t \vv = \nu\Delta \vv - (\vv \cdot \vn) \vv  -\vn q+ \vF,
\end{equation*}
with $q=q_1+q_2+q_{3}$ and $\vF=\vF_0+\nu\vF_1+\vF_2+\vF_3$.\\

Now that we have obtained the expressions defining $q$ and $\vF$, we must prove the size estimates on the pressure and the force. This will be done in the following two lemmas.
%%%%%%%%%%%%%%%%%%%%%%%%%%%%%%%%%%%%%%%%%%%%%%
\begin{Lemme}\label{Lemme3232pression}The pressure $q$ is a function such that $q\in L^{3/2}_tL^{3/2}_x(\Omega_0)$.
\end{Lemme}
%%%%%%%%%%%%%%%%%%%%%%%%%%%%%%%%%%%%%%%%%%%%%%
\textbf{\textit{Proof.}} We recall that $q=q_{1}+q_{2}+q_{3}$, where the expressions $q_{1}$, $q_{2}$ and $q_{3}$ were given in (\ref{Definition_q1}), (\ref{definition_q2}) and (\ref{definition_q3}) respectively, \emph{i.e.}:
$$q_{1}=-\frac{1}{\Delta}(\Div(\vor \wedge \psi\vu )), \qquad q_{2}=\frac{1}{\Delta} \Div( \psi \vec{A})\quad  \mbox{ and } \quad q_{3}=-\frac{1}{2}|\vv|^{2}.$$
We will study each one of these terms separately. The first and the last term (\emph{i.e.} $q_1$ and $q_3$) are very easy to deal with. 
\begin{itemize}
\item For $q_{1}$,   we just write the following estimates :
\begin{eqnarray*}
\| q_{1}\|_{L^{3/2}(B_{x_{0}, \rho_0})} &=&\left\|\frac{1}{\Delta}(\vn \cdot ( \vor \wedge \psi\vu ))\right\|_{L^{3/2}(B_{x_{0}, \rho_0})}\leq  C_{\rho_0}\left\|\frac{1}{\Delta} \vn \cdot ( \vor \wedge \psi \vu)\right\|_{L^{9/5}( \mathbb{R}^3)}\\
&\leq & C_{\rho_0} \|   \vor \wedge \psi \vu\|_{L^{9/8}(\mathbb{R}^{3})}\leq  C_{\rho_0}\|\vor\|_{L^{2}(B_{x_{0}, \rho})} \|\psi \vu\|_{L^{18/7}(B_{x_{0}, \rho})} \\
&\leq & C_{\rho_0} \|\vn\wedge \vu\|_{L^{2}(B_{x_{0}, \rho})} \|\psi\|_{L^{\infty}}   \|\vu\|_{L^{18/7}(B_{x_{0}, \rho})} \\
&\leq & C_{\rho_0}\|\psi\|_{L^{\infty}}\|\vu\|^{4/3}_{H^{1}(B_{x_{0}, \rho})} \|\vu\|^{2/3}_{L^{2}(B_{x_{0}, \rho})},
\end{eqnarray*}
where we used the H\"older inequality, Hardy-Littlewood-Sobolev inequalities, the support properties of $\psi$ and an interpolation estimate. Now, integrating with respect to the time variable we find
$$\| q_{1}\|_{L^{3/2}_tL^{3/2}_x(\Omega_0)}\leq C_{\rho_0, \psi}\|\vu\|^{4/3}_{L^{2}(I_0, H^{1}(B_{x_{0},\rho}))}  \|\vu\|^{2/3}_{L^{\infty}(I_0,L^{2}(B_{x_{0},\rho}))}<+\infty.$$
\item  For $q_{3}=-\frac{1}{2}|\vv|^{2}$, we just write
\begin{eqnarray*}
\| q_{3}\|_{L^{3/2}(B_{x_{0}, \rho_0})} =  C\|\vv\|^{2}_{L^{3}(B_{x_{0}, \rho_0})} \leq C \|\vv\|_{H^{1}(B_{x_{0}, \rho_0})} \|\vv\|_{L^{2}(B_{x_{0}, \rho_0})},
\end{eqnarray*}
and taking the $L^{3/2}$-norm in time we have
$$ \| q_{3}\|_{L^{3/2}_tL^{3/2}_x(\Omega_0)}\leq  C_{\rho_{0}}   \|\vv\|_{L^2(I_0,H^{1}(B_{x_{0}, \rho_0}))} \|\vv\|_{L^\infty(I_0,L^{2}(B_{x_{0}, \rho_0}))} <+\infty.$$\\[2mm]
\end{itemize}
Thus, the main term we have to study is the term $q_2$. Recall that  $q_{2}=\frac{1}{\Delta} \Div(\psi \vec{A})$ with  
$$\vec A= (\vn \wedge \vv) \wedge\vec \beta +  (\vn\wedge\vec\beta) \wedge \vv + (\vn\wedge\vec\beta) \wedge \vec\beta$$ 
where $\vec \beta$ was defined in (\ref{Definition_Beta}). Thus, applying the H\"older inequality and the Hardy-Littlewood-Sobolev inequality we write
$$\|q_{2}\|_{L^{3/2}(B_{x_{0}, \rho_0})}\leq  C_{\rho_0} \|q_{2}\|_{L^{9/5}(B_{x_{0}, \rho_0})}\leq C_{\rho_0}\|\psi\vec A\|_{L^{9/8}(B_{x_{0}, \rho_{0}})}.$$
Moreover,
\begin{eqnarray*} 
\|\psi\vec A\|_{L^{9/8}(B_{x_{0}, \rho_{0}})}&\leq&   \|\vn \wedge \vv \|_{L^2(B_{x_{0}, \rho_{0}})} \| \psi \vec \beta \|_{L^{18/7}(B_{x_{0}, \rho_{0}})}\\
&&+ \quad  \|\vn\wedge\vec\beta\|_{L^2(B_{x_{0}, \rho_{0}})} \|\psi \vv \|_{L^{18/7}(B_{x_{0}, \rho_{0}})}\\
&&+\quad \|\vn\wedge\vec\beta\|_{L^2(B_{x_{0},\rho_{0}})}  \|\psi \vec\beta\|_{L^{18/7}(B_{x_{0},\rho_{0}})}.
\end{eqnarray*}
Now, using an interpolation inequality we have
\begin{eqnarray}
\|\psi\vec A\|_{L^{9/8}(B_{x_{0}, \rho_{0}})}&\leq& \|\vn \wedge \vv \|_{L^2(B_{x_{0},\rho_{0}})}  \| \psi \vec \beta \|_{L^{2}(B_{x_{0},\rho_{0}})}^{2/3}  \| \psi \vec \beta \|_{H^1(B_{x_{0},\rho_{0}})}^{1/3}\nonumber \\ 
& &+ \quad \|\vn \wedge \vec\beta \|_{L^2(B_{x_{0},\rho_{0}})}  \| \psi \vv \|_{L^{2}(B_{x_{0},\rho_{0}})}^{2/3}  \| \psi \vv \|_{H^1(B_{x_{0},\rho_{0}})}^{1/3}\label{Estimation_A}\\ 
&&+\quad   \|\vn \wedge  \vec\beta \|_{L^2(B_{x_{0},\rho_{0}})}  \| \psi \vec \beta \|_{L^{2}(B_{x_{0},\rho_{0}})}^{2/3}  \| \psi \vec \beta \|_{H^1(B_{x_{0},\rho_{0}})}^{1/3}.\nonumber
\end{eqnarray}
Recalling that $\vec{\beta}=-\vw$ and applying Corollary \ref{Proprietes_DifferenceUV} we obtain 
\begin{eqnarray*}
\|\psi\vec \beta\|_{L^2(B_{x_{0},\rho_{0}})}&\leq &C_{\rho, \psi}\|\vu\|_{L^2(B_{x_{0},\rho_{0}})},
\end{eqnarray*}
and in a similar manner we have
$$\|\vn\otimes\vec\beta\|_{L^2(B_{x_{0},\rho_{0}})}\leq C_{\rho, \psi}  \|\vu\|_{H^1(B_{x_{0},\rho_{0}})}.$$
Thus, with these previous inequalities, we obtain for the expression (\ref{Estimation_A}) the following estimate
\begin{eqnarray*} 
\|q_{2}\|_{L^{3/2}(B_{x_{0}, \rho_0})}&\leq &C_{\rho_{0}}\|\psi\vec A\|_{L^{9/8}(B_{x_{0},\rho_{0}})}\\
&\leq& C_{\rho_{0},\psi}   \|\vv \|_{H^1(B_{x_{0},\rho_{0}})}  \|\vu\|_{L^2(B_{x_{0},\rho_{0}})}^{2/3} \|\vu \|_{H^1(B_{x_{0},\rho_{0}})}^{1/3}  \\ 
&&+ C_{\rho_{0},\psi} \|\vu\|_{H^1(B_{x_{0},\rho_{0}})} \|\vv\|_{L^2(B_{x_{0},\rho_{0}})}^{2/3} \|\vv \|_{H^1(B_{x_{0},\rho_{0}})}^{1/3}  \\ 
&&+  C_{\rho_{0},\psi}   \|\vu\|_{H^1(B_{x_{0},\rho_{0}})}^{4/3}  \|\vu\|_{L^2(B_{x_{0},\rho_{0}})}^{2/3},
\end{eqnarray*}
and finally, since we have $\vv\in L^{\infty}_{t}L^{2}_{x}\cap L^{2}_{t}H^{1}_{x}(\Omega_{0})$ and $\vu\in L^{\infty}_{t}L^{2}_{x}\cap L^{2}_{t}H^{1}_{x}(\Omega)$, we obtain $q_2\in L^{3/2}_{t}L^{3/2}_{x}(\Omega_0)$.\\
 
We have thus proven that $q=q_{1}+q_{2}+q_{3}\in L^{3/2}_{t}L^{3/2}_{x}(\Omega_0)$ and the proof of Lemma \ref{Lemme3232pression} is finished. \hfill $\blacksquare$\\

%%%%%%%%%%%%%%%%%%%%%%%%%%%%%%%%%%%%%%%%%%%%%%
\begin{Lemme}\label{LemmeL2L2Force}The force $\vF$ is such that $\Div(\vF)=0$ and $\vF-\vec F_0\in L^2_tL^2_x( \Omega_0)$.
\end{Lemme}
%%%%%%%%%%%%%%%%%%%%%%%%%%%%%%%%%%%%%%%%%%%%%%
\textbf{ \textit{Proof.}} We already saw that $\Div(\vF)=0$. Thus, the only point to check is that $\vF-\vec F_0\in L^2_tL^2_x( \Omega_0)$.  Again, we will study each term of $\vF-\vF_0=\nu\vF_1+\vF_2+\vF_3$ separately. 
\begin{itemize}
\item For $\vF_{1}$, we start by recalling its definition given in (\ref{Definition_F1}):
\begin{eqnarray*}
\vF_1 &= &-\frac{1}{\Delta} \vn \wedge\left[\vn \wedge(\Delta\psi \vu) - \vn (\Delta \psi)\wedge \vu\right] + 2 \frac{1}{\Delta} \vn \wedge\left[ \sum_{j=1}^3\partial_{x_j}\left(  \vn \wedge\left( (\partial_{x_j} \psi)  \vu\right)\right)\right] \\
& &- 2\frac{1}{\Delta} \vn \wedge \left[\sum_{j=1}^3\partial_{x_j}\left(\vn (\partial_{x_j}\psi) \wedge \vu\right)\right].
\end{eqnarray*}
We have then
\begin{eqnarray*} \|\vec F_1\|_{L^2(B_{x_{0}, \rho_0})}&\leq& C\bigg(\|\Delta\psi\, \vu\|_{L^{2}(B_{x_{0}, \rho_0})}+  \|\vn (\Delta \psi)\wedge \vu\|_{L^{6/5}(B_{x_{0}, \rho})}\\
& &+\sum_{j=1}^3\|\vn\wedge((\partial_j\psi)\vu)\|_{L^{2}(B_{x_{0}, \rho})}+\sum_{j=1}^3\|  (\vn\partial_j\psi)\wedge \vu\|_{L^{2}(B_{x_{0}, \rho})}\bigg)\\
\end{eqnarray*}
Thus, with the properties of $\psi$ we obtain
$$ \|\vF_1\|_{L^2_tL^2_x(\Omega_0)}\leq C_\psi \|\vu\|_{L^2(I_0,H^1(B_{x_{0},\rho})}\leq C_\psi \|\vu\|_{L^2_t H^1_x(\Omega)}<+\infty.$$
%%%%%%%%%%%%%%%%%%%%%%%%%%%%%%%%%%%%%%%%%%%%%%
\item For $\vF_{2}$ we have by the formula (\ref{Definition_F2}):
$$\|\vF_2\|_{L^{2}(B_{x_{0}, \rho_0})}=\left\|-\frac{1}{\Delta} \vn \wedge \left((\vn \psi)\wedge (\vor \wedge \vu)\right)\right\|_{L^{2}(B_{x_{0}, \rho_0})}.$$
Here, we are going to apply the same arguments used in the study of the quantity $\vec{w}$ treated in (\ref{Formule_Kernel}). Indeed, we can write
$$\left\|-\frac{1}{\Delta} \vn \wedge \left((\vn \psi)\wedge (\vor \wedge \vu)\right)\right\|_{L^{2}(B_{x_{0}, \rho_0})}\leq C_{\rho_0}\left\|\frac{1}{\Delta} \vn \left((\vn \psi)\cdot(\vor\wedge\vu)\right)\right\|_{L^{\infty}(B_{x_{0},\rho_0)}},$$
but since, for $x\in B(x_{0},\rho_0)$, we have
$$\left|\frac{1}{\Delta} \vn \left((\vn \psi)\cdot(\vor\wedge\vu)\right)(t,x)\right| \leq C\int_{\{\rho_1-\rho_0<\vert x-y\vert, \, y\in B_{x_{0},\rho}\}}\frac{|\vn \psi(t,y)|}{|x-y|^{2}} |\vor(t,y)| |\vu(t,y)|dy,$$
we obtain the following bound
$$\left|\frac{1}{\Delta} \vn \left((\vn \psi)\cdot(\vor\wedge\vu)\right)(t,x)\right|\leq C_{\rho_0}\|\vn \psi\|_{L^{\infty}(B_{x_{0},\rho})}\int_{B_{x_{0},\rho}}|\vor(t,y)| |\vu(t,y)|dy.$$
Now,  taking the $L^{2}$ norm in the time variable we have,
\begin{eqnarray*}
\| \vF_2 \|_{L^{2}_{t}L^{2}_{x}(\Omega_0)}&\leq &C_{\rho_0,\psi}\left\|\vor\right\|_{L^{2}(I_0,L^{2}(B_{x_{0}, \rho}))}\left\|\vu\right\|_{L^{\infty}(I_0,L^{2}(B_{x_{0}, \rho}))}\\
& \leq & C_{\rho_0,\psi}\left\|\vu\right\|_{L^{2}_{t}H^{1}_{x}(\Omega)}\left\|\vu\right\|_{L^{\infty}_{t}L^{2}_{x}(\Omega)}<+\infty.
\end{eqnarray*}
%%%%%%%%%%%%%%%%%%%%%%%%%%%%%%%%%%%%%%%%%%%%%%
\item For $\vF_{3}$ we have $\vF_{3}=-\psi \vec{A}+\vn \frac{1}{\Delta} \Div(\psi \vec{A})$. Recalling that $\psi(t,x)=\phi(t)\Phi(x)$, we introduce a new cut-off function $\eta\in\mathcal{D}(\mathbb{R}^3)$ which is equal to $1$ on $B_{x_{0}, \rho_{3}}$ and is supported within $B_{x_{0}, \rho_{2}}$ with $\rho_0<\rho_3<\rho_2<\rho_1<\rho$, so that on $\Omega_0$ we have
$$\vF_{3}= -\phi\eta \vec{A}+\vn \frac{1}{\Delta} \Div(\phi\eta \vec{A})+\vn \frac{1}{\Delta} \Div\big(\psi(1-\eta) \vec{A}\big),$$
thus we can write
\begin{equation}\label{Estimation_F3}
\|\vF_3\|_{L^2_tL^2_x(\Omega_0)} \leq \|\phi\eta \vec{A}\|_{L^2_t L^2_x(\Omega_0)} +\left\|\vn \frac{1}{\Delta} \Div(\phi\eta \vec{A})\right\|_{L^2_t L^2_x(\Omega_0)}+\left\|\vn \frac{1}{\Delta} \Div\big(\psi(1-\eta) \vec{A}\big)\right\|_{L^2_tL^2_x(\Omega_0)}.
\end{equation}
We study the first term $\|\phi\eta \vec{A}\|_{L^2_t L^2_x(\Omega_0)}$ and with the definition of $\vec{A}$ given in (\ref{Definition_A}) we have the following estimates in the space variable
\begin{eqnarray}
\|\phi\eta \vec{A}\|_{L^2(B_{x_{0}, \rho_{0}})}&\leq &\big\|\phi\eta \big((\vn \wedge \vv)\wedge \vec{\beta}\big)\big\|_{L^2(B_{x_{0}, \rho_{0}})}+\big\|\phi\eta\big((\vn \wedge \vec{\beta})\wedge \vv\big)\big\|_{L^2(B_{x_{0}, \rho_{0}})}\label{MajorationA1}\\
& &+\big\|\phi\eta \big((\vn \wedge \vec{\beta})\wedge \vec{\beta}\big)\big\|_{L^2(B_{x_{0}, \rho_{0}})}\nonumber\\
&\leq & \|\vn \wedge \vv\|_{L^2(B_{x_{0}, \rho_{0}})}\|\phi\eta\vec{\beta} \|_{L^\infty(B_{x_{0}, \rho_{0}})}+\|\vv\|_{L^2(B_{x_{0}, \rho_{0}})}\|\phi\eta(\vn \wedge\vec{\beta} )\|_{L^\infty(B_{x_{0}, \rho_{0}})}\nonumber\\
& & + \|\vn \wedge\vec{\beta}\|_{L^2(B_{x_{0}, \rho_{0}})}\|\phi\eta\vec{\beta}\|_{L^\infty(B_{x_{0}, \rho_{0}})}.\nonumber
\end{eqnarray} 
Thus, since $\vec{\beta}=-\vw$ and since by Proposition \ref{Proposition_1} $\vw\in L^\infty_t {\rm Lip}_x(\Omega_{0})$, we have the inequalities
\begin{eqnarray*}
\|\phi\eta\vec{\beta} \|_{L^\infty(B_{x_{0}, \rho_{0}})}&\leq &\|\vec{\beta} \|_{L^\infty(B_{x_{0}, \rho_{0}})}\leq C_{\rho_{0}}\|\vn \psi\|_{L^{2}(B_{x_{0},\rho})} \|\vu\|_{L^{2}(B_{x_{0},\rho})},\\ 
\|\phi\eta(\vn \wedge\vec{\beta} )\|_{L^\infty(B_{x_{0}, \rho_{0}})}&\leq &\|\vn \wedge\vec{\beta}\|_{L^\infty(B_{x_{0}, \rho_{0}})}\leq C_{\rho_{0}}\|\vn \psi\|_{L^{2}(B_{x_{0},\rho})} \|\vu\|_{L^{2}(B_{x_{0},\rho})},\\
\|\vn \wedge\vec{\beta}\|_{L^2(B_{x_{0}, \rho_{0}})}&\leq &C_{\rho_{0}}\|\vn \psi\|_{L^{2}(B_{x_{0},\rho})} \|\vu\|_{L^{2}(B_{x_{0},\rho})},
\end{eqnarray*}
after an integration in the time variable we obtain
\begin{eqnarray} 
\|\phi\eta\vec{A}\|_{L^2_t L^2_x(\Omega_{0})}&\leq& C_{\rho, \psi}   \|\vn \wedge \vv \|_{L^2_t L^2_x(\Omega_{0})} \|\vu\|_{L^\infty_t L^2_x(\Omega)}\nonumber\\
&& + C_{\rho, \psi} \|\vv \|_{L^2_t L^2_x(\Omega_{0})} \|\vu\|_{L^\infty_t L^2_x(\Omega)}+C_{\rho, \psi}\|\vu\|_{L^2_{t} L^2_x(\Omega)}\|\vu\|_{L^\infty_t L^2_x(\Omega)}\nonumber\\ 
&\leq&   C_\rho \|\vu\|_{L^\infty_t L^2_x(\Omega)} \left(\|\vu\|_{L^2_t H^1_x(\Omega)} +\|\vv\|_{L^2_t H^1_x(\Omega)}\right)<+\infty.\label{MajorationA2} 
\end{eqnarray}

Now we study the second term of (\ref{Estimation_F3}) and by the support properties of $\eta$ we can write
\begin{eqnarray*}
\left\|\vn \frac{1}{\Delta} \Div(\phi\eta \vec{A})\right\|_{L^2(B_{x_{0}, \rho_{0}})}&\leq & \left\|\vn \frac{1}{\Delta} \Div(\phi\eta \vec{A})\right\|_{L^2(\mathbb{R}^{3})}\leq C\|\phi\eta \vec{A}\|_{L^2(\mathbb{R}^{3})}\\
&\leq & C\|\phi\eta \vec{A}\|_{L^2(B_{x_{0}, \rho_{2}})}.
\end{eqnarray*}
Since $\rho_{2}<\rho_{1}<\rho$ and thus $B_{x_{0}, \rho_{2}}\subset B_{x_{0}, \rho}$, we can apply the same arguments used in (\ref{MajorationA1})-(\ref{MajorationA2}) to obtain the following estimate (see also Remark \ref{Remarque3}):
\begin{equation}\label{MajorationA3}
\left\|\vn \frac{1}{\Delta} \Div(\phi\eta \vec{A})\right\|_{L^{2}_{t}L^2_{x}(\Omega_{0})}\leq C_\rho \|\vu\|_{L^\infty_t L^2_x(\Omega)} \left(\|\vu\|_{L^2_t H^1_x(\Omega)} +\|\vv\|_{L^2_t H^1_x(\Omega)}\right)<+\infty.
\end{equation}
%%%%%%%%%%%%%%%%%%%%%%%%%%%%%%%%%%%%%%%%%%%%%%
Due to the support properties of $\psi(1-\eta) $, the last term of (\ref{Estimation_F3}) can not be treated in the same manner as before and we have for $(t,x)\in \Omega_0$
\begin{equation*}\begin{split}  
\left| \vn \frac{1}{\Delta} \Div(\psi(1-\eta)  \vec{A})(t,x) \right|\leq& C \int_{\{\rho_3-\rho_0<\vert x-y\vert,\, y\in B_{x_{0},\rho}\}} \frac{1}{\vert x-y\vert^3} \vert\psi(t,y)\vert\,  \vert\vec A(t,y)\vert\, dy  \\ 
\leq C \|\psi\|_{L^\infty(B_{x_{0},\rho})} \frac{1}{(\rho_3-\rho_0)^3} &\left( \|\vn\wedge\vv\|_{L^2(B_{x_{0},\rho})} \|\psi\vec\beta\|_{L^2(B_{x_{0},\rho})}+\|\vn\wedge \vec\beta\|_{L^2(B_{x_{0},\rho})} \|\psi\vv\|_{L^2(B_{x_{0},\rho})}\right.\\ 
& \quad + \left. \|\vn\wedge\vec\beta\|_{L^2(B_{x_{0},\rho})} \|\psi\vec\beta\|_{L^2(B_{x_{0},\rho})}\right).
\end{split}
\end{equation*}
Now, since $\vec{\beta}=-\vw$, we apply Corollary \ref{Proprietes_DifferenceUV} to obtain:
\begin{equation*}
\begin{split}
\left| \vn \frac{1}{\Delta} \Div(\psi(1-\eta)  \vec{A})(t,x) \right|\leq C_{\rho,\psi} &\left(  \|\vv\|_{H^1(B_{x_{0},\rho})} \|\vu\|_{L^2(B_{x_{0},\rho})} +\|\vu\|_{H^1(B_{x_{0},\rho})} \|\vv\|_{L^2(B_{x_{0},\rho})} \right.\\ 
&\quad \left.+\|\vu\|_{H^1(B_{x_{0},\rho})} \|\vu\|_{L^2(B_{x_{0},\rho})} \right),
\end{split}
\end{equation*}
so that
\begin{eqnarray}
\left\|\vn \frac{1}{\Delta} \Div(\psi(1-\eta) \vec{A})\right\|_{L^2_tL^2_x(\Omega_0)}&\leq& C_{\rho, \psi} \bigg(\|\vv\|_{L^2_t H^1_x(\Omega)} \|\vu\|_{L^\infty_t L^2_x(\Omega)} + \|\vu\|_{L^2_t H^1_x(\Omega)} \|\vv\|_{L^\infty_t L^2_x(\Omega)}\nonumber\\
&& +\|\vu\|_{L^2_t H^1_x(\Omega)} \|\vu\|_{L^\infty_t L^2_x(\Omega)}\bigg)<+\infty.\label{MajorationA4}
\end{eqnarray}
With the inequalities (\ref{MajorationA2}), (\ref{MajorationA3}) and (\ref{MajorationA4}) we finally obtain that $\vF_{3}\in L^2_tL^2_x(\Omega_0)$.
\end{itemize}
Gathering together the estimates over $\vF_1$, $\vF_2$ and $\vF_3$ we finally have that 
$$\|\vF-\vF_{0}\|_{L^2_tL^2_x(\Omega_0)}<+\infty,$$
and the proof of Lemma \ref{LemmeL2L2Force} is finished. \hfill $\blacksquare$\\

At this point, with Lemmas \ref{Lemme3232pression} and \ref{LemmeL2L2Force} we have proven that the new variable 
$$\vv=-\frac{1}{\Delta} \vn \wedge (\psi \vn \wedge \vu),$$ 
satisfies \emph{locally} the following Navier--Stokes system
$$\partial_t \vv =\Delta \vv - (\vv \cdot \vn)\vv - \vn q+\vF,$$
where $q$ is a pressure such that $q\in L^{3/2}_tL^{3/2}_x(\Omega_0)$ and  $\vF$ is a force such that $\Div(\vF)=0$ and $\vF-\vF_0\in L^2_tL^2_x(\Omega_0)$, which is the conclusion of Proposition \ref{Proposition_2}.\hfill $\blacksquare$\\

It might be interesting to notice that, if $\vec u$ is more regular, we have a better estimate on $\vec F-\vec F_0$ :
%%%%%%%%%%%%%%%%%%%%%%%%%%%%%%%%%%%%%%%%%%%%%%
\begin{Lemme}\label{LemmeL2L2Force_bis} 
If moreover $\vu\in L^\infty_t H^1_x(\Omega)\cap  L^2_t H^2_x(\Omega)$, then $\vF-\vec F_0\in L^2_tH^1_x( \Omega_0)$.
\end{Lemme}
%%%%%%%%%%%%%%%%%%%%%%%%%%%%%%%%%%%%%%%%%%%%%%
\textbf{\textit{Proof.}} Same proof as for  Lemma \ref{LemmeL2L2Force}.  \hfill $\blacksquare$
%%%%%%%%%%%%%%%%%%%%%%%%%%%%%%%%%%%%%%%%%%%%%%
\subsection{The case $\vf\in L^{10/7}_t L^{10/7}_x$}\label{Secc_Existence_Dissipativite}
In this section, we  shall prove Proposition \ref{Proposition_40} and point $1)$ in Theorem \ref{theo_CLM}, which we recall now :
%%%%%%%%%%%%%%%%%%%%%%%%%%%%%%%%%%%%%%%%%%%%%%
\begin{Proposition}\label{Proposition_401}Let  $\Omega$ be a bounded subset of $\mathbb{R}\times \mathbb{R}^{3}$ of the form (\ref{Set_Definition}) and assume that $\vu\in L^{\infty}_{t}L^{2}_{x}(\Omega)\cap L_{t}^{2}H_{x}^{1}(\Omega)$ with $\Div(\vu)=0$ and $p\in \mathcal{D}'(\Omega)$ are solutions of the Navier--Stokes equations (\ref{Equation_NS}) on $\Omega$ where $\vf\in L^{10/7}_t L^{10/7}_x(\Omega)$ and $\Div(\vf)=0$.\\
\begin{itemize}
\item[1)]  Let $\gamma\in \mathcal{D}(\mathbb{R})$ and  $\theta \in \mathcal{D}(\mathbb{R}^{3})$ be two smooth functions such that $\displaystyle{\int_{\mathbb{R}}} \gamma(t)dt=\displaystyle{\int_{\mathbb{R}^{3}}} \theta(x)dx=1$, $supp(\gamma)\subset ]-1,1[$ and $supp(\theta)\subset B(0,1)$. We set for $\alpha, \varepsilon >0$ the functions $\gamma_{\alpha}(t)=\frac{1}{\alpha}\gamma(\frac{t}{\alpha})$ and $\theta_{\varepsilon}=\frac{1}{\varepsilon^{3}}\theta(\frac{x}{\varepsilon  })$ and we define $\varphi_{\alpha, \varepsilon}(t,x)=\gamma_{\alpha}(t)\theta_{\varepsilon}(x)$. Then, if the cylinder  $Q_{t_0,x_0,r}$ is contained in $ \Omega$, the distributions $\vu \ast \varphi_{\alpha, \varepsilon}$ and $p\ast \varphi_{\alpha, \varepsilon}$ are well defined in the set $Q_{t_0,x_0,r/4}\subset \Omega$ for $0<\alpha<r_0^2/2 $ and $0<\varepsilon<r_0/2$. Moreover, the limit
$$\underset{\varepsilon \to 0}{\lim}\,\underset{\alpha \to 0}{\lim}\, \Div\left[\big(p\ast \varphi_{\alpha, \varepsilon}\big) \times \big(\vu \ast \varphi_{\alpha, \varepsilon}\big)\right],$$
exists in $\mathcal{D}'$ and does not depend on the choice of the functions $\gamma$ and $\theta$. \\

\item[2)] Let $$\underset{\varepsilon \to 0}{\lim}\,\underset{\alpha \to 0}{\lim}\, \Div\left[\big(p\ast \varphi_{\alpha, \varepsilon}\big) \times \big(\vu \ast \varphi_{\alpha, \varepsilon}\big)\right]= \langle\Div(p\vu)\rangle.$$ Then  the quantity
\begin{equation}\label{Definition_M}
\begin{split}
M=-\partial_{t}|\vu|^{2}+\nu \Delta |\vu|^{2}-2\nu|\vn \otimes \vu|^{2}-&\Div(|\vu|^{2}\vu)\\& -2\langle \Div(p\vu)\rangle+ 2 \vu\cdot \vf,
\end{split}
\end{equation}
 is well-defined as a distribution.
\end{itemize}
\end{Proposition}
%%%%%%%%%%%%%%%%%%%%%%%%%%%%%%%%%%%%%%%%%%%%%%
We make here two remarks. The first one is about the hypothesis $\vf \in  L^{10/7}_t L^{10/7}_x$:  this is enough to give a sense to the product $\vu\cdot \vf$ in (\ref{Definition_M}) since we have $\vu \in L^{10/3}_t L^{10/3}_x$. Our second remark concerns the proof of this proposition:  we will need the properties of the companion equation studied in the previous section in order to obtain the existence of these limits and the crucial link between the properties of $\vu$ and $\vv$ is given in Lemma \ref{main_lemma} below.\\

\textit{\textbf{Proof.}} 
\begin{itemize}
\item[1)] We start with the first point of the proposition and for this we denote by $\vu_{\alpha, \varepsilon}$ the function defined by $\vu_{\alpha, \varepsilon}=\vu \ast \varphi_{\alpha, \varepsilon}$. Since this is a regular function in the time and space variables we can write
$$ \partial_{t} |\vu_{\alpha, \varepsilon}|^{2}=  2 \vu_{\alpha, \varepsilon} \cdot \partial_{t}\vu_{\alpha, \varepsilon}$$ so that
\begin{eqnarray*}
\partial_{t} |\vu_{\alpha, \varepsilon}|^{2}&=& 2 \vu_{\alpha, \varepsilon}\cdot  \big(\partial_{t}\vu \ast  \varphi_{\alpha, \varepsilon}\big)= 2 \vu_{\alpha, \varepsilon} \cdot  \left(\nu\Delta \vu  -(\vu \cdot \vn )\vu -\vn p+\vf\right)\ast \varphi_{\alpha, \varepsilon}\\
&=&2 \vu_{\alpha, \varepsilon} \cdot \nu\Delta \vu_{\alpha, \varepsilon} - 2 \vu_{\alpha, \varepsilon} \cdot \left([(\vu \cdot \vn )\vu]\ast \varphi_{\alpha, \varepsilon} \right)- 2 \vu_{\alpha, \varepsilon} \cdot (\vn p \ast \varphi_{\alpha, \varepsilon}) + 2 \vu_{\alpha, \varepsilon} \cdot (\vf \ast \varphi_{\alpha, \varepsilon}).
\end{eqnarray*}
Using the fact that $\Div(\vu)=0$, we have:
\begin{equation}\label{Dissipative_1}\begin{split}
\partial_{t} |\vu_{\alpha, \varepsilon}|^{2}=&\nu \Delta |\vu_{\alpha, \varepsilon}|^{2}-2\nu |\vn \otimes \vu_{\alpha, \varepsilon}|^{2}-2\vu_{\alpha, \varepsilon}\cdot( \left[ \Div(\vu\otimes \vu) \right]\ast \varphi_{\alpha, \varepsilon})\\&-2\, \Div\left[(p\ast \varphi_{\alpha, \varepsilon}) \vu_{\alpha, \varepsilon} \right]+2 \vu_{\alpha, \varepsilon}\cdot (\vf \ast \varphi_{\alpha, \varepsilon}).\end{split}
\end{equation}
We shall now take $\Omega_0$ in the previous section large enough to contain $Q_{t_0,x_0,r/2}$ and we consider the Navier--Stokes companion equation (\ref{Equation_V0}) on $\vv$ with pressure $q$ and force $\vF$.  Let us remark that since $\vf$ belongs to $L^{10/7}_t L^{10/7}_x(\Omega)$, by Lemma \ref{LemmeForce0} we have  $\vF_{0}\in L^{10/7}_t L^{10/7}_x(\Omega_{0})$, moreover since by Proposition \ref{Proposition_2} we have $\vF-\vF_{0} \in L^{2}_t L^{2}_x(\Omega_{0})$ and we obtain that $\vF\in L^{10/7}_t L^{10/7}_x(\Omega_{0})$.\\

Now, if we denote by $\vv_{\alpha, \varepsilon}$ the function given by $\vv_{\alpha, \varepsilon}=\vv \ast \varphi_{\alpha, \varepsilon}$, by the same arguments used above, we obtain the following equation
\begin{equation}\label{Dissipative_2}\begin{split}
\partial_{t} |\vv_{\alpha, \varepsilon}|^{2}=&\nu \Delta |\vv_{\alpha, \varepsilon}|^{2}-2\nu |\vn \otimes \vv_{\alpha, \varepsilon}|^{2}-2\vv_{\alpha, \varepsilon}\cdot (\left[ \Div(\vv\otimes \vv)\right] \ast \varphi_{\alpha, \varepsilon})\\ &-2\, \Div\left[(q\ast \varphi_{\alpha, \varepsilon}) \vv_{\alpha, \varepsilon} \right]+2 \vv_{\alpha, \varepsilon}\cdot (\vF \ast \varphi_{\alpha, \varepsilon}).
\end{split}\end{equation}
We aim to study the convergence of expressions (\ref{Dissipative_1}) and (\ref{Dissipative_2})  as the parameters $\alpha$ and $\varepsilon$ vanish and we will use the properties of $\vv$ to deduce the limits on $\vu$. However, these limits must be treated very carefully and we will first make $\alpha\longrightarrow 0$ and only then we will take the limit $\varepsilon\longrightarrow 0$. \\

To begin our study, we remark that it is easy to deal with the convergence of some of the terms contained in those formulas. For the sake of simplicity, we will adopt the following notation: $\vu_{\varepsilon}=\vu\ast \theta_{\varepsilon}$ and we will denote by $\bar{Q}$ the space-time cylinder $Q_{t_0,x_0,r/4}\subset \Omega_{0}$.
%%%%%%%%%%%%%%%%%%%%%%%%%%%%%%%%%%%%%%%%%%%%%%
\begin{Lemme}\label{LemDissipative1} 
We have the following strong convergence
\begin{itemize}
\item[1)] $\vu_{\alpha, \varepsilon}\underset{\alpha\to 0}{\longrightarrow} \vu_{\varepsilon}$ in $L^{2}_{t}L^{2}_{x}(\bar{Q})$ and in $L^{2}_{t}\dot{H}^{1}_{x}(\bar{Q})$,
\item[2)]  $\vn \otimes \vu_{\alpha, \varepsilon}\underset{\alpha\to 0}{\longrightarrow} \vn \otimes \vu_{\varepsilon}$  in $L^{2}_{t}L^{2}_{x}(\bar{Q})$,
\item[3)] $(\vu \otimes \vu)\ast \varphi_{\alpha, \varepsilon}\underset{\alpha \to 0}{\longrightarrow}  (\vu \otimes \vu)\ast \theta_{\varepsilon}$  in $L^{2}_{t}L^{2}_{x}(\bar{Q})$,
\item[4)] $\vf \ast \varphi_{\alpha, \varepsilon}\underset{\alpha \to 0}{\longrightarrow} \vf\ast \theta_{\varepsilon}$  in $L^{10/7}_t L^{10/7}_x(\bar{Q})$.
 \end{itemize}
\end{Lemme}
%%%%%%%%%%%%%%%%%%%%%%%%%%%%%%%%%%%%%%%%%%%%%%
Due to the properties of the functions $\vu$, $\vf$ and with the definition of the function $\varphi _{\alpha, \varepsilon}=\gamma_{\alpha} \theta_{\varepsilon}$, the proof of this lemma  is straightforward.
\begin{Remarque}\label{Remarque5} Conclusions of Lemma  \ref{LemDissipative1} can be obtained in the same manner for the new variable $\vv$: indeed, as shown in Propositions \ref{Proposition_1} and \ref{Proposition_2},  $\vv$ and $\vF$ have a similar behavior as $\vu$ and $\vf$, see also Lemma \ref{LemmeForce0}. However, we have another convergence for $\vv$ that will be very useful in the sequel: $q\ast \varphi_{\alpha,\epsilon}$ is strongly convergent when $\alpha \longrightarrow 0$ to $q\ast \theta_\epsilon$ in $L^{3/2}_t L^{3/2}_{x}(\bar Q)$.
 \end{Remarque}
Now we can pass to the limit $\alpha \longrightarrow 0$ for all the terms in equality (\ref{Dissipative_1}) except for the term involving $p$. But the limit for this term must exist, as (due to the equality)  it is equal to a sum of terms that do have a limit.  Thus, we obtain
\begin{equation}\label{DefLimite11}
\begin{split}
\partial_{t} |\vu_{\varepsilon}|^{2}=\nu \Delta |\vu_{\varepsilon}|^{2}-2\nu |\vn \otimes \vu_{\varepsilon}|^{2}-2\vu_{\varepsilon}\cdot( \left[ \Div(\vu\otimes \vu)\right]  \ast \theta_{\varepsilon})\\- 2\underset{\alpha\to 0}{\lim}\ \Div\left[(p\ast \varphi_{\alpha, \varepsilon}) \vu_{\alpha, \varepsilon} \right]+2 \vu_{\varepsilon}\cdot( \vf \ast \theta_{\varepsilon}),
\end{split}
\end{equation}
and, in the same manner we have for the function $\vv$:
\begin{equation}\label{DefLimite12}
\begin{split}
 \partial_{t} |\vv_{\varepsilon}|^{2}=\nu \Delta |\vv_{\varepsilon}|^{2}-2\nu |\vn \otimes \vv_{\varepsilon}|^{2}-2\vv_{\varepsilon}\cdot (\left[ \Div(\vv\otimes \vv)\right] \ast \theta_{\varepsilon})\\ -2 \, \Div\left[(q\ast \theta_{\varepsilon}) \vv_{ \varepsilon} \right]+2 \vv_{\varepsilon}\cdot ( \vF \ast \theta_{\varepsilon}).
 \end{split}
\end{equation}
Note in particular that since we have the stronger convergence $\underset{\alpha \to 0}{\lim}q\ast \varphi_{\alpha,\epsilon}=q\ast \theta_\epsilon$ in $L^{3/2}_t L^{3/2}_{x}$ we can write 
$\Div\left[(q\ast \theta_{\varepsilon}) \vv_{ \varepsilon} \right]$ for the couple $(\vv, q)$ in (\ref{DefLimite12}) instead of $\underset{\alpha\to 0}{\lim}\ \Div\left[(p\ast \varphi_{\alpha, \varepsilon}) \vu_{\alpha, \varepsilon} \right]$ for $(\vu, p)$ in (\ref{DefLimite11}).\\

At this point we define two quantities that will help us to pass to the limit $\varepsilon \to 0$:
\begin{eqnarray}
\mu_{\varepsilon}&=& 2 \vu_{\varepsilon}\cdot( \left[ \Div(\vu\otimes \vu)\right]  \ast \theta_{\varepsilon}) -\Div(|\vu|^{2}\vu)\label{Defaut1}\\[3mm]
\eta_{\varepsilon}&=& 2 \vv_{\varepsilon}\cdot( \left[ \Div(\vv\otimes \vv)\right]  \ast \theta_{\varepsilon}) -\ \Div(|\vv|^{2}\vv),\label{Defaut2}  
\end{eqnarray}
we will see with these quantities how to link the behavior of $\vv$ to the behavior of $\vu$.  But before this, we will need the following lemma which states some strong convergence results in the space variable. 
%%%%%%%%%%%%%%%%%%%%%%%%%%%%%%%%%%%%%%%%%%%%%%
\begin{Lemme}\label{LemDissipative4} We have the following strong convergence
\begin{itemize}
\item[1)]   $\vu_{\varepsilon}\underset{\varepsilon\to 0}{\longrightarrow} \vu$ in $L^{2}_{t}L^{2}_{x}(\bar{Q})$  and in $L^{2}_{t}\dot{H}^{1}_{x}(\bar{Q})$,
\item[2)]
 $\vn \otimes \vu_{\varepsilon}\underset{\varepsilon \to 0}{\longrightarrow}\vn \otimes \vu$  in $L^{2}_{t}L^{2}_{x}(\bar{Q})$,
\item[3)]    $\vf \ast\theta_{\varepsilon}\underset{\varepsilon \to 0}{\longrightarrow} \vf$  in $L^{10/7}_t L^{10/7}_x(\bar{Q})$.
\end{itemize}
\end{Lemme}
%%%%%%%%%%%%%%%%%%%%%%%%%%%%%%%%%%%%%%%%%%%%%%
Again, the proof of this lemma is straightforward. Of course, conclusions of Lemma  \ref{LemDissipative4} can be obtained in the same manner for the new variable $\vv$. But, as pointed out in Remark \ref{Remarque5}, we have another convergence for $\vv$, in the space variable this time:  $q\ast \theta_{\epsilon}$ is strongly convergent to $q$ in $L^{3/2}_t L^{3/2}_{x}(\bar Q)$.\\
 
With the help of this lemma, passing to the limit $\varepsilon \longrightarrow 0$ in (\ref{DefLimite11}) and (\ref{DefLimite12}) we then have:
\begin{equation}\label{DefLimite13}  
\begin{split}
\partial_{t} |\vu|^{2}=&\nu \Delta |\vu|^{2}-2\nu |\vn \otimes \vu|^{2}-  \Div(|\vu|^{2}\vu) + 2 \vu \cdot \vf \\
&- \underset{\varepsilon\to 0}{\lim} \left( \mu_{\varepsilon}  +2 \underset{\alpha\to 0}{\lim}\, \Div\left[(p\ast \varphi_{\alpha, \varepsilon})\times \vu_{\alpha, \varepsilon} \right] \right),
\end{split}
\end{equation}
and 
\begin{equation}\label{DefLimite14}  
\begin{split}
\partial_{t} |\vv|^{2}=&\nu \Delta |\vv|^{2}-2\nu |\vn \otimes \vv|^{2}-  \Div(|\vv|^{2}\vv) +2 \vv \cdot \vF\\ 
& - \underset{\varepsilon\to 0}{\lim}\, \eta_{\varepsilon}  -2\,   \Div(q \vv).
\end{split}
\end{equation}

Although we have obtained at this stage similar equations for $\vu$ and $\vv$ the situation is quite different: we do not have any information about the pressure $p\in \mathcal{D}'(\Omega)$ but we do have a much better behavior for the new pressure $q$, since  $q\in L^{3/2}_{t}L^{3/2}_{x}(\bar{Q})$ (as it  was proved in Proposition \ref{Proposition_2}).\\

The end of the proof now relies on the following lemma that relates the behavior of $\mu_\epsilon$ to the behavior of $\eta_\epsilon$ :
\begin{Lemme}[Key lemma]\label{main_lemma}
For $\mu_{\varepsilon}$ and $\nu_{\varepsilon}$ defined in (\ref{Defaut1}) and (\ref{Defaut2}) respectively, we have the following convergence in $\mathcal{D}'(\bar Q)$:
$$ \underset{\varepsilon\to 0}{\lim}\, \eta_{\varepsilon}-\mu_\epsilon=0.$$ 
\end{Lemme}
Before going into the details of the proof, let us explain the consequences of this lemma. If we have $ \underset{\varepsilon\to 0}{\lim}\, \eta_{\varepsilon}-\mu_\epsilon=0$, then the existence of $ \underset{\varepsilon\to 0}{\lim}\, \eta_{\varepsilon}$ will imply the existence of $ \underset{\varepsilon\to 0}{\lim}\, \mu_{\varepsilon}$ and we can give a sense to the quantity
$$\underset{\varepsilon \to 0}{\lim}\,\underset{\alpha \to 0}{\lim}\, \Div\left[\big(p\ast \varphi_{\alpha, \varepsilon}\big) \times \big(\vu \ast \varphi_{\alpha, \varepsilon}\big)\right],$$
as all the remaining terms of (\ref{DefLimite13}) do exists since $\vu \in L^\infty_{t}L^2_{x}(\bar Q)\cap L^2_{t}H^1_{x}(\bar Q)$ and $\vf\in L^{10/7}_t L^{10/7}_x(\bar Q)$ and thus the first point of Proposition \ref{Proposition_401} will be proven. \\

But the existence of $\underset{\varepsilon\to 0}{\lim}\, \eta_{\varepsilon}$ is given by the properties of the companion equation: indeed, all the terms of identity (\ref{DefLimite14}) exists since $\vv \in L^\infty_{t}L^2_{x}(\bar Q)\cap L^2_{t}H^1_{x}(\bar Q)$, $\vF\in L^{10/7}_t L^{10/7}_x(\bar Q)$ with $q\in L^{3/2}_{t}L^{3/2}_{x}(\bar Q)$ and thus the limit $\underset{\varepsilon\to 0}{\lim}\, \eta_{\varepsilon}$ exists. See also Remark \ref{Remarque0} and the references \cite{Vasseur, LEMA2}.\\

As we can see, Lemma \ref{main_lemma} explains how to link the properties of the companion equation to the original Navier--Stokes equation. Once the scope of this lemma is clear, we turn our attention to its proof which relies on an idea of Duchon \& Robert \cite{DR1}:
%%%%%%%%%%%%%%%%%%%%%%%%%%%%%%%%%%%%%%%%%%%%%%
\begin{Proposition}\label{theo:duchon13} Let  $\Omega$ be a bounded subset of $\mathbb{R}\times \mathbb{R}^{3}$ of the form (\ref{Set_Definition}) and let $\vec u\in L^\infty_{t}L^2_{x}(\Omega)\cap L^2_{t}H^1_{x}(\Omega)$, with $ \Div(\vu)=0.$ Let $\theta \in \mathcal{D}(\mathbb{R}^{3})$ be a smooth function  such that $ \displaystyle{\int_{\mathbb{R}^{3}}} \theta(x)dx=1$  and $supp(\theta)\subset B(0,1)$. We set for $  \varepsilon >0$ the function $\theta_{\varepsilon}=\frac{1}{\varepsilon^{3}}\theta(\frac{x}{\varepsilon  })$. Then, if the cylinder  $Q_{t_0,x_0,r}$ is contained in $ \Omega$,  we define the following  distributions on $Q_{t_0,x_0,r/4}\subset \Omega$ for   $0<\varepsilon<r_0/2$ :
\begin{eqnarray*}
\mu_\epsilon&=&2 (\vu\ast\theta_{\varepsilon})\cdot( \left[ \Div(\vu\otimes \vu)\right]  \ast \theta_{\varepsilon}) -\Div(|\vu|^{2}\vu)\\
R_\varepsilon &= &\sum_{k=1}^3 \int_{\mathbb{R}^{3}} \partial_k\theta_\varepsilon(y)  \big(u_k(x-y)-u_k(x)\big) \vert\vec u(x-y)-\vec u(x)\vert^2\, dy\\
S_\varepsilon&=& \sum_{k=1}^3 \int_{\mathbb{R}^{3}} \partial_k\theta_\varepsilon(y)  \big(u_k(x-y)-u_k(x)\big) \big(\vec u(x-y)-\vec u(x))\cdot(\vec u_\epsilon(x)-\vec u(x)\big)\, dy.
\end{eqnarray*}
Then, we have the following limit
$$\underset{\varepsilon\to 0}{\lim}\, \mu_{\varepsilon}+R_\varepsilon-2S_\varepsilon=0.$$
\end{Proposition}   
%%%%%%%%%%%%%%%%%%%%%%%%%%%%%%%%%%%%%%%%%%%%%%
The proof of this proposition is given in the appendix. \\

Once we have this result at hand, we will use it in the following way: we introduce the notation 
$$\tau_{z}[\vec H](t,x)=\vec  H(t,x-z)- \vec H(t,x),$$ 
where $\vec H:\mathbb{R}\times \mathbb{R}^{3}\longrightarrow \mathbb{R}^{3}$ is a function and $z\in \mathbb{R}^{3}$ is a vector. Then we define
\begin{eqnarray*}
T_{\varepsilon}(\vec U, \vec V, \vec W)&=& -\sum_{k=1}^{3}\int_{\mathbb{R}^{3}} \partial_{x_{k}} \theta_{\varepsilon}(y) \tau_{y}[U_{k}](t,x)\big(\tau_{y}[\vec V](t,x)\cdot\tau_{y}[\vec W](t,x)\big)dy\\
& & + 2\sum_{k=1}^{3}\int_{\mathbb{R}^{3}} \partial_{x_{k}} \theta_{\varepsilon}(y) \tau_{y}[U_{k}](t,x)\bigg(\tau_{y}[\vec V](t,x)\cdot \big(\theta_{\varepsilon}\ast \vec W (t,x)-\vec W(t,x)\big)  \bigg)dy,
\end{eqnarray*}
and we remark that we have $R_\varepsilon-2S_\varepsilon=-T_{\varepsilon}(\vec U,\vec U,\vec U)$, observe moreover that the operator $T_{\varepsilon}$ is a trilinear operator.\\

Thus, from Proposition \ref{theo:duchon13}, we see that the solutions $\vu$ and $\vv$ of the Navier--Stokes equations on $\Omega_0$ that have been discussed previously satisfy 
$$ \underset{\varepsilon\to 0}{\lim}\, \mu_{\varepsilon}-\eta_\varepsilon+T_{\varepsilon}(\vv,\vv,\vv)-T_{\varepsilon}(\vu,\vu,\vu)=0,$$
with $\mu_\varepsilon$ and $\eta_\varepsilon$ defined by (\ref{Defaut1}) and (\ref{Defaut2}). Then proving Lemma \ref{main_lemma} amounts to prove that
$$ \underset{\varepsilon\to 0}{\lim}\,  T_{\varepsilon}(\vv,\vv,\vv)-T_{\varepsilon}(\vu,\vu,\vu)=0.$$
For this we write
\begin{equation}\label{Decomposition_Tri}
T_{\varepsilon}(\vv, \vv, \vv)-T_{\varepsilon}(\vu, \vu, \vu)=T_{\varepsilon}(\vec v-\vec u, \vv, \vv)+T_{\varepsilon}(\vu, \vec v-\vec u, \vv)+T_{\varepsilon}(\vu, \vu,  \vec v-\vec u),
\end{equation}
and then we conclude with the following lemma:
%%%%%%%%%%%%%%%%%%%%%%%%%%%%%%%%%%%%%%%%%%%%%%
\begin{Lemme}\label{trilinear}If $\vec U$, $\vec V$ and $\vec W$ belong to $L^3_t L^3_x(\Omega_0)$ and if at least one of them belong to $L^\infty_t {\rm Lip}_x(\Omega_0)$, then we have 
$$\lim_{\varepsilon\rightarrow 0}  T_{\varepsilon}(\vec U, \vec V, \vec W)=0,$$
in $L^1_t L^1_x(Q_{t_0,x_0,r/4})$.
\end{Lemme}
%%%%%%%%%%%%%%%%%%%%%%%%%%%%%%%%%%%%%%%%%%%%%%
\textit{\textbf{Proof.}}  Recalling that $\theta_{\varepsilon}=\frac{1}{\varepsilon^{3}}\theta(\frac{x}{\varepsilon})$ and $supp(\theta_{\varepsilon})\subset B(0,\varepsilon)$ we can write
\begin{eqnarray*}
\left|T_{\varepsilon}(\vec U, \vec V, \vec W)(t,x)\right|&\leq& C\frac{1}{\varepsilon^4}  \int_{\{\vert y\vert<\varepsilon\}} \vert \vec U(t,x)-\vec U(t,x-y)\vert\,   \vert \vec V(t,x)-\vec V(t,x-y)\vert\\
&& \times  \vert \vec W(t,x)-\vec W(t,x-y)\vert\, dy  \\ 
&+& \frac{C}{\varepsilon^7}  \int_{\{\vert y\vert<\varepsilon\}} \vert \vec U(t,x)-\vec U(t,x-y)\vert\,   \vert \vec V(t,x)-\vec V(t,x-y)\vert\,  dy\\
& &\times  \int_{\{\vert z\vert<\varepsilon\}} \vert \vec W(t,x)-\vec W(t,x-z)\vert\, dz, 
\end{eqnarray*}
and applying the H\"older inequality we have
\begin{eqnarray}
\left|T_{\varepsilon}(\vec U, \vec V, \vec W)(t,x)\right|&\leq &\frac{C}{\varepsilon^4} \left(\int_{\{|y|<\varepsilon\}} \!\!\!\! \vert \vec U(t,x)-\vec U(t,x-y)\vert^3\, dy\right)^{\frac{1}{3}}\nonumber\\
& & \times  \left(\int_{\{|y|<\varepsilon\}} \!\!\!\!  \vert \vec V(t,x)-\vec  V(t,x-y)\vert^3\, dy\right)^{\frac{1}{3}}\label{OperateurT}\\
& &\times \left(\int_{\{|y|<\varepsilon\}} \!\!\!\!  \vert \vec W(t,x)-\vec W(t,x-y)\vert^3\, dy\right)^{\frac{1}{3}}.\nonumber
\end{eqnarray}
Now, if $\vec U\in L^3_t L^3_x(\Omega_0)$, then we may write
$$\frac{1}{\varepsilon^3} \iint_{Q_{t_0,x_0,r/4}} \int_{\{|y|<\varepsilon\}} \vert \vec U(t,x)-\vec U(t,x-y)\vert^3\, dy\, dt\, dx=\int_{\{|z|<1\}}\iint_{Q_{t_0,x_0,r/4}}  \vert \vec U(t,x)-\vec U(t,x-\varepsilon z)\vert^3 \, dt\, dx\, dz. 
$$
and we can check easily (by dominated convergence) that the right-hand term goes to $0$ as $\varepsilon\rightarrow 0$.  Moreover if we have $\vec U\in L^\infty_t {\rm Lip}_x(\Omega_0)$, then we obtain
  $$ \frac{1}{\varepsilon^6} \iint_{Q_{t_0,x_0,r/4}} \int_{\{|y|<\varepsilon\}}  \vert \vec U(t,x)-\vec U(t,x-y)\vert^3\, dy\, dt\, dx \leq  C \|\vn\otimes\vec U\|_{L_{t}^{\infty}L_{x}^{\infty}}^3  \vert Q_{t_0,x_0,r/4}\vert.
$$
With these inequalities, if $\vec U\in L^\infty_t {\rm Lip}_x(\Omega_0)$ and if $\vec V\in L^3_t L^3_x(\Omega_0)$, $\vec W \in L^3_t L^3_x(\Omega_0)$ (say), integrating (\ref{OperateurT}) over $Q_{t_0,x_0,r/4}$ we obtain
\begin{eqnarray*}
\iint_{Q_{t_0,x_0,r/4}}\left|T_{\varepsilon}(\vec U, \vec V, \vec W)(t,x)\right|dtdx &\leq & C\|\vn\otimes\vec U\|_{L_{t}^{\infty}L_{x}^{\infty}} \vert Q_{t_0,x_0,r/4}\vert^{\frac{1}{3}}\\
& & \times \left(\int_{\{|z|<1\}}\iint_{Q_{t_0,x_0,r/4}}  \vert \vec V(t,x)-\vec V(t,x-\varepsilon z)\vert^3 \, dt\, dx\, dz\right)^{\frac{1}{3}}\\
& & \times \left(\int_{\{|z|<1\}}\iint_{Q_{t_0,x_0,r/4}}  \vert \vec W(t,x)-\vec W(t,x-\varepsilon z)\vert^3 \, dt\, dx\, dz\right)^{\frac{1}{3}},
\end{eqnarray*}
we have then 
$$\lim_{\varepsilon\rightarrow 0} \iint_{Q_{t_0,x_0,r/4}}\left|T_{\varepsilon}(\vec U, \vec V, \vec W)(t,x)\right|dtdx =0,$$
and the lemma is proven.    \hfill $\blacksquare$\\
 
Let us finish the proof of  Lemma \ref{main_lemma}. Since over $\Omega_{0}$, we have $\vu, \vv\in L^\infty_t L^2_x(\Omega_0)\cap L^2_t H^1_x(\Omega_0)$ we easily obtain that $\vu, \vv\in L^3_t L^3_x(\Omega_0)$, moreover we have by Proposition \ref{Proposition_1} that $\vv-\vu\in L^\infty_t {\rm Lip}_x(\Omega_0)$ and we can apply Lemma \ref{trilinear} to each term of (\ref{Decomposition_Tri}), thus using Proposition \ref{theo:duchon13} we obtain 
$$\underset{\varepsilon\to 0}{\lim}\, \eta_{\varepsilon}-\mu_\epsilon=0,$$ 
and Lemma \ref{main_lemma} is proven.\hfill $\blacksquare$\\

As said before, we already know that $\underset{\varepsilon\to 0}{\lim}\, \eta_{\varepsilon}$, thus  $\underset{\varepsilon\to 0}{\lim}\, \mu_{\varepsilon}$ exists. As we know that 
$$\underset{\varepsilon\to 0}{\lim} \left( \mu_{\varepsilon}  +2 \underset{\alpha\to 0}{\lim}\, \Div\left[(p\ast \varphi_{\alpha, \varepsilon})\times \vu_{\alpha, \varepsilon}\right] \right),$$ 
exists, we find that $ \underset{\varepsilon\to 0}{\lim}   \underset{\alpha\to 0}{\lim}\, \Div\left[(p\ast \varphi_{\alpha, \varepsilon})\times \vu_{\alpha, \varepsilon} \right]$ exists. Moreover since $\underset{\varepsilon\to 0}{\lim}\, \eta_{\varepsilon}=\underset{\varepsilon\to 0}{\lim}\,\mu_\epsilon$ and using the identities (\ref{DefLimite13}) and (\ref{DefLimite14}) we have
\begin{eqnarray*}
2\underset{\varepsilon\to 0}{\lim}   \underset{\alpha\to 0}{\lim}\, \Div\left[(p\ast \varphi_{\alpha, \varepsilon})\times \vu_{\alpha, \varepsilon} \right]&=& -\partial_{t} |\vu|^{2}+\nu \Delta |\vu|^{2}-2\nu |\vn \otimes \vu|^{2}\\&&-  \Div(|\vu|^{2}\vu) + 2 \vu \cdot \vf   - \underset{\varepsilon\to 0}{\lim}  \mu_{\varepsilon} \\&=& 2\, \Div(q\vv)+\partial_t(|\vv|^2-|\vu|^2)+\nu\Delta(|\vu|^2-|\vv|^2)\\ && +2\nu(|\vn\otimes\vv|^2-|\vn\otimes\vu|^2) +2 (\vu \cdot\vf-\vv \cdot\vF)\\&& +\Div( |\vv|^2\vv-|\vu|^2\vu).
\end{eqnarray*}
Clearly, this limit does not depend on the choice of the functions $\gamma$ and $\theta$. We have obtained that the limit
$$\underset{\varepsilon \to 0}{\lim}\,\underset{\alpha \to 0}{\lim}\, \Div\left[\big(p\ast \varphi_{\alpha, \varepsilon}\big) \times \big(\vu \ast \varphi_{\alpha, \varepsilon}\big)\right],$$
exists in $\mathcal{D}'$ which is the conclusion of the first point of the proposition.

\item[2)] The second point of Proposition \ref{Proposition_401} is easy to prove. Indeed, since we have
\begin{eqnarray*}
2\underset{\varepsilon\to 0}{\lim}   \underset{\alpha\to 0}{\lim}\, \Div\left[(p\ast \varphi_{\alpha, \varepsilon}) \times \vu_{\alpha, \varepsilon} \right]&=& -\partial_{t} |\vu|^{2}+\nu \Delta |\vu|^{2}-2\nu |\vn \otimes \vu|^{2}\\&&-  \Div(|\vu|^{2}\vu) + 2 \vu \cdot \vf   - \underset{\varepsilon\to 0}{\lim}  \mu_{\varepsilon},
\end{eqnarray*}
we obtain that the distribution $M$ introduced in (\ref{Definition_M}) is in fact given by the following identity
$$M=\underset{\varepsilon\to 0}{\lim}\, \mu_{\varepsilon}.$$
\end{itemize}
The proposition is proven.\hfill $\blacksquare$
\begin{Remarque}\label{Remarque44}Since $\vf\in L^{10/7}_{t}L^{10/7}_{x}$,  if we assume $\vu\in L^{4}_{t}L^{4}_{x}(\Omega)$, then it is easy to see that $M=0$. Indeed, in this setting we have $\vv \in L^{4}_{t}L^{4}_{x}(\Omega_{0})$ and if we consider the companion equation (\ref{Equation_V0}) on $\vv$ we obtain a Navier--Stokes equation with a pressure $q\in L^{3/2}_{t}L^{3/2}_{x}(\Omega_{0})$ and a force $\vF\in L^{10/7}_t L^{10/7}_x(\Omega_{0})$, then following \cite{LEMA2} we have $\underset{\varepsilon\to 0}{\lim}\, \eta_{\varepsilon}=0$, but by Lemma \ref{main_lemma} we have $\underset{\varepsilon\to 0}{\lim}\, \eta_{\varepsilon}=\underset{\varepsilon\to 0}{\lim}\, \mu_{\varepsilon}$ and thus we obtain $M=0$.
\end{Remarque}

%%%%%%%%%%%%%%%%%%%%%%%%%%%%%%%%%%%%%%%%%%%%%%
%%%%%%%%%%%%%%%%%%%%%%%%%%%%%%%%%%%%%%%%%%%%%%
\subsection{The case $\vf\in L^{2}_{t}L^{2}_{x}$}
In this section, we shall prove point $2)$ in Theorem \ref{theo_CLM}, which we recall now :

%%%%%%%%%%%%%%%%%%%%%%%%%%%%%%%%%%%%%%%%%%%%%%
\begin{Proposition}\label{Proposition_402}Let  $\Omega$ be a bounded subset of $\mathbb{R}\times \mathbb{R}^{3}$ of the form (\ref{Set_Definition}) and assume that $\vu\in L^{\infty}_{t}L^{2}_{x}(\Omega)\cap L_{t}^{2}H_{x}^{1}(\Omega)$ with $\Div(\vu)=0$ and $p\in \mathcal{D}'(\Omega)$ are solutions of the Navier--Stokes equations on  $\Omega$ :
\begin{equation*}
\partial_{t} \vu=\nu\Delta \vu -(\vu \cdot \vn )\vu -\vn p+\vf,
\end{equation*}
where  $\Div(\vf)=0$. Assume that:
\begin{itemize}
\item $\vf\in L^2_tL^2_x(\Omega)$,
\item $\vu$ is dissipative: the distribution $M$ defined by 
\begin{equation*}
\begin{split}
M=-\partial_{t}|\vu|^{2}+\nu \Delta |\vu|^{2}-&2\nu|\vn \otimes \vu|^{2}- \Div(|\vu|^{2}\vu)\\& -2\langle \Div(p\vu)\rangle + 2 \vu\cdot \vf\end{split}
\end{equation*}
is a non-negative locally finite Borel measure on $\Omega$.
\end{itemize}
There exists positive constants $\epsilon^*>0$ and $\tau_1>5$  which depend  only on $\nu$  such that, if $(t_0,x_0)\in \Omega$ and 
$$\limsup_{r\rightarrow 0} \frac{1}{r}\iint_{]t_0-r^2,t_0+r^2[\times B_{x_0,r}} \vert\vec\nabla\otimes \vec u(s,y)\vert^2\, ds\, dy<\epsilon^*,$$ 
 then  there exists a    neighborhood $\mathcal{Q}=Q_{t_{0}, x_{0}, \bar{r}}$  of $(t_0,x_0)$  such that $\mathds{1}_{\mathcal{Q}} \,\vec u\in \mathcal{M}^{3,\tau_1}_2$. 
\end{Proposition}
%%%%%%%%%%%%%%%%%%%%%%%%%%%%%%%%%%%%%%%%%%%%%%
\textit{\textbf{Proof.}} We shall again take $\Omega_0$  as in the previous sections,  large enough to contain $(t_0,x_0)$ and we consider the Navier--Stokes companion equation (\ref{Equation_V0}) on $\vv$ with pressure $q$ and force $\vF$. \\

Let us first remark that, since we have now the hypothesis $\vf\in L^2_t L^2_x$,  we obtain by Lemma \ref{LemmeForce0} that $\vF_{0}$ given in (\ref{Definition_F0}) belongs to the space $L^2_t L^2_x(\Omega_0)$ and as by Proposition \ref{Proposition_2} we have $\vF-\vF_{0}\in L^2_t L^2_x(\Omega_0)$, we obtain that the global force $\vF$ belongs to $L^2_t L^2_x(\Omega_0)$. Then just as in the previous section, we will exploit the properties of $\vv$ to deduce the wished results on $\vu$. Indeed we have:
\begin{itemize}
\item $q\in L^{3/2}_{t}L^{3/2}_{x}(\Omega_0)$,
\item ${\bf 1}_{\Omega_0} \vF   \in L^2_t L^2_x\subset \mathcal{M}^{10/7,2}_2$,
\item $\vv$ is suitable, 
\item there exists a positive constant $\epsilon^*>0$ which depend only on $\nu$ such that, if $(t_0,x_0)\in \Omega_0$ we have  
\begin{equation}\label{PetitGrad}
\limsup_{r\rightarrow 0} \frac{1}{r}\iint_{]t_0-r^2,t_0+r^2[\times B_{x_0,r}} \vert\vec\nabla\otimes \vec v(s,y)\vert^2\, ds\, dy<\epsilon^*.
\end{equation}
\end{itemize}

The two first points are straightforward, for the suitability of $\vv$, let us recall that the distribution $M$ associated to $\vu$ has been computed in the previous sections as $M=\underset{\varepsilon\rightarrow 0}{\lim}\mu_\varepsilon$ and by the dissipativity hypothesis 
we have $M\geq 0$. As we have seen that 
$\underset{\varepsilon\rightarrow 0}{\lim} \mu_\varepsilon=\underset{\varepsilon\rightarrow 0}{\lim}\eta_\varepsilon$, we find that  $\underset{\varepsilon\rightarrow 0}{\lim} \eta_\varepsilon$ is a non-negative locally finite Borel measure on $\Omega_0$ and thus $\vv$ is suitable.\\

The last point, \emph{i.e.} the condition (\ref{PetitGrad}), can be deduced from the hypothesis on $\vu$ stated above: we use the fact that $\vw=\vv-\vu$ belongs to $L^\infty_t {\rm Lip}_x(\Omega_0)$ and then from the inequality
\begin{equation*}\begin{split}
\left| \|  \mathds{1}_{Q_{t_{0}, x_{0}, r}}\vn\otimes\vv\|_{L^{2}_{t}L^{2}_{x}} - \|\mathds{1}_{Q_{t_{0}, x_{0}, r}}\vn\otimes\vu\|_{L^{2}_{t}L^{2}_{x}}\right| \leq& \| \mathds{1}_{Q_{t_{0}, x_{0}, r}}\vn\otimes\vw\|_{L^{2}_{t}L^{2}_{x}}\\ 
\leq& C r^{5/2} \|\vn\otimes\vw\|_{L^\infty_tL^\infty_x(Q_{t_{0}, x_{0}, r} )}, 
 \end{split}\end{equation*}
we find that
\begin{equation*}
\begin{split} 
\limsup_{r\rightarrow 0} \frac{1}{r}\iint_{]t_0-r^2,t_0+r^2[\times B_{x_0,r}}& \vert\vec\nabla\otimes \vec u(s,y)\vert^2\, ds\, dy =\\ &\limsup_{r\rightarrow 0} \frac{1}{r}\iint_{]t_0-r^2,t_0+r^2[\times B_{x_0,r}} \vert\vec\nabla\otimes \vec v(s,y)\vert^2\, ds\, dy,
\end{split}
\end{equation*}
thus, if the first limit (on $\vu$) is less than $\epsilon^*$, the second limit (on $\vv$) is still less than $\epsilon^*$.\\

As we can see, we gathered enough information on the solution $\vv$ of the companion equation in order to apply point $2)$ in Kukavica's theorem (Theorem \ref{theo_kukavica}): there exists $\tau_{1}>5$ and a small neighborhood $\mathcal{Q}=Q_{t_{0}, x_{0}, \bar{r}}\subset \Omega_{0}$  of $(t_0,x_0)$ such that $\mathds{1}_{\mathcal{Q}} \,\vec v\in \mathcal{M}^{3,\tau_1}_2$ and $\mathds{1}_{\mathcal{Q}}\,  q\in \mathcal{M}^{3/2,\tau_1/2}_2$. \\
 
We have obtained interesting properties for the solution $\vv$ of the companion equation, but now we must go back to $\vu$: we already have $\mathds{1}_{\mathcal{Q}} \,\vec v\in \mathcal{M}^{3,\tau_1}_2$ in a small  neighborhood $\mathcal{Q}=Q_{t_{0}, x_{0}, \bar{r}}$  of $(t_0,x_0)$. 
As $\vw\in L^\infty_{t}L^\infty_{x}(\Omega_0)\subset L^{\tau_1}_{t}L^{\tau_1}_{x}(\Omega_0)$, we have $\mathds{1}_{\mathcal{Q}} \,\vec  w\in \mathcal{M}^{\tau_{1},\tau_1}_2$  and since $\tau_{1}>3$ by Remark \ref{RemarqueMorrey} we obtain that $\mathds{1}_{\mathcal{Q}} \,\vec  w\in \mathcal{M}^{3,\tau_1}_2$ as well. Thus $\mathds{1}_{\mathcal{Q}} \,\vec u\in \mathcal{M}^{3,\tau_1}_2$ and the proposition is proven.  \hfill $\blacksquare$
%%%%%%%%%%%%%%%%%%%%%%%%%%%%%%%%%%%%%%%%%%%%%%
%%%%%%%%%%%%%%%%%%%%%%%%%%%%%%%%%%%%%%%%%%%%%%
\subsection{The case $\vf\in L_{t}^{ 2} H_{x}^1$}

In this section, we shall prove point $3)$ in Theorem \ref{theo_CLM} with the help of the following proposition.

%%%%%%%%%%%%%%%%%%%%%%%%%%%%%%%%%%%%%%%%%%%%%%
\begin{Proposition}\label{Proposition_403} Let  $\Omega$ be a bounded subset of $\mathbb{R}\times \mathbb{R}^{3}$ of the form (\ref{Set_Definition}) and assume that $\vu\in L^{\infty}_{t}L^{2}_{x}(\Omega)\cap L_{t}^{2}\dot{H}_{x}^{1}(\Omega)$ with $\Div(\vu)=0$ and $p\in \mathcal{D}'(\Omega)$ are solutions of the Navier--Stokes equations on  $\Omega$ :
\begin{equation*}
\partial_{t} \vu=\nu\Delta \vu -(\vu \cdot \vn )\vu -\vn p+\vf,
\end{equation*}
where  $\Div(\vf)=0$. Assume that   there exists a neighborhood $\mathcal{Q}=Q_{t_{0}, x_{0}, \bar{r}}$  of $(t_0,x_0)\in \Omega$  such that 
\begin{itemize}
\item  $\mathds{1}_{\mathcal{Q}}\, \vec u\in \mathcal{M}^{3,\tau_1}_2$ for some $\tau_1>5$,
\item  $\mathds{1}_{\mathcal{Q}} \, \vf \in L^2_tH^1_x$.  
\end{itemize}
Then there exists $r'<\bar{r}$  such that $\vu$ is bounded on $Q_{t_0,x_0,r'}$. In particular,  the point $(t_0,x_0)$ is regular.
\end{Proposition}
%%%%%%%%%%%%%%%%%%%%%%%%%%%%%%%%%%%%%%%%%%%%%%
\textit{\textbf{Proof.}}
Just apply O'Leary's theorem (Theorem \ref{Theoreme_2}).\hfill $\blacksquare$\\

With this last proposition, we have finished the proof the Theorem \ref{theo_CLM} from which we deduce Theorem \ref{Theorem_1}.

\begin{Remarque} As we don't have any information over the pressure $p$, we can't just apply the last point of Kukavica's theorem and we need to invoke O'Leary's result.
\end{Remarque}
This last remark says that we have in fact a better result for the companion equation:
%%%%%%%%%%%%%%%%%%%%%%%%%%%%%%%%%%%%%%%%%%%%%%
\begin{Proposition}\label{Proposition_404} Let  $\Omega$ be a bounded subset of $\mathbb{R}\times \mathbb{R}^{3}$ of the form (\ref{Set_Definition}) and assume that $\vu\in L^{\infty}_{t}L^{2}_{x}(\Omega)\cap L_{t}^{2}\dot{H}_{x}^{1}(\Omega)$ with $\Div(\vu)=0$ and $p\in \mathcal{D}'(\Omega)$ are solutions of the Navier--Stokes equations on  $\Omega$ :
\begin{equation*}
\partial_{t} \vu=\nu\Delta \vu -(\vu \cdot \vn )\vu -\vn p+\vf,
\end{equation*}
where  $\Div(\vf)=0$. Assume that  there exists a neighborhood $\mathcal{Q}=Q_{t_{0}, x_{0}, \bar{r}}$ of $(t_0,x_0)\in \Omega$ such that $\mathds{1}_{\mathcal{Q}}\, \vec u\in \mathcal{M}^{3,\tau_1}_2$ for some $\tau_1>5$ and $\mathds{1}_{\mathcal{Q}} \, \vf \in L^2_tH^1_x$.  Assume that for the companion equation on $\vv$, $q$ and $\vF$ we have as well $\mathds{1}_{\mathcal{Q}}\,  q\in \mathcal{M}^{3/2,\tau_1/2}_2$. Then there exists $r'<\bar{r}$  such that $\vv$ is H\"olderian (with parabolic H\"older regularity exponent $\eta>0$) on $Q_{t_0,x_0,r'}$. 
\end{Proposition}
%%%%%%%%%%%%%%%%%%%%%%%%%%%%%%%%%%%%%%%%%%%%%%
\textit{\textbf{Proof.}} First, we apply Proposition \ref{Proposition_403} and find a cylinder  $Q_{t_0,x_0,r_1}$ on which $\vu$ is bounded and thus by the Serrin regularity theory, on a smaller cylinder $Q_{t_0,x_0,r_2}$ we have that $\vu\in L^\infty_t H^1_{x}\cap L^2_t H^2_{x}$. Next, since $\mathds{1}_{\mathcal{Q}} \, \vf \in L^2_tH^1_x$ by Lemma \ref{LemmeForce0} we obtain that $\vF_0\in L^2_{t} H^1_{x}(Q_{t_0,x_0,r_2})$, thus with these informations for $\vu$ and $\vf$, we consider the companion equation and applying the same computations as for Proposition \ref{Proposition_2} (see also Lemma \ref{LemmeL2L2Force_bis}), we then see that, on a still smaller cylinder $Q_{t_0,x_0,r_3}$ we have $\vF-\vF_0\in L^2_{t} H^1_{x}$. This gives for the global force $\vF\in L^2_{t} H^1_{x}(Q_{t_0,x_0,r_3})$  from which we deduce that $\mathds{1}_{Q_{t_0,x_0,r_3}} \vF\in \mathcal{M}^{10/7,10/3}_2$.  As $10/3>5/2$, we may now apply point $3)$ in Kukavica's theorem (Theorem \ref{theo_kukavica}) and conclude. \hfill $\blacksquare$\\

As a conclusion, we can see that the obstruction to time regularity lies inside the harmonic correction $\vw$, a fact that is of course reminiscent of Serrin's counter-example.

%%%%%%%%%%%%%%%%%%%%%%%%%%%%%%%%%%%%%%%%%%%%%%
\appendix
\section*{Appendix}
\textit{\textbf{Proof of Proposition \ref{theo:duchon13}.}} We start the proof by the remark that, as $\Div(\vu)=0$, we have on $Q_{t_0,x_0,r/4}$ the identity
$$\sum_{k=1}^3 \int_{\mathbb{R}^{3}} \partial_k\theta_\epsilon(y)  (u_k(x-y)-u_k(x))\, dy=0,$$
so that we can rewrite $R_\epsilon$ in the following manner
\begin{eqnarray*}
R_\epsilon&=&\sum_{k=1}^3 \int_{\mathbb{R}^{3}} \partial_k\theta_\epsilon(y)  \left(u_k(x-y)-u_k(x)\right) \left(\vert\vec u(x-y)\vert^2-2 \vec u(x-y)\cdot \vec u(x)\right)\, dy\\
&=& \sum_{k=1}^3 \int_{\mathbb{R}^{3}} \partial_k\theta_\epsilon(y)  \left(u_k(x-y)-u_k(x)\right) \vert\vec u(x-y)\vert^2\, dy\\
& & -2\sum_{k=1}^3 \int_{\mathbb{R}^{3}} \partial_k\theta_\epsilon(y)  \left(u_k(x-y)-u_k(x)\right) \vec u(x-y)\cdot \vec u(x)\,dy\\
R_\epsilon&=& \Div\left( \theta_\epsilon* (\vert\vec u\vert^2\vec u)-(\theta_\epsilon*\vert\vec u\vert^2)\vec u\right)-2\sum_{k=1}^3 \int_{\mathbb{R}^{3}} \partial_k\theta_\epsilon(y)  \left(u_k(x-y)-u_k(x)\right) \vec u(x-y)\cdot \vec u(x)\,dy.
\end{eqnarray*}
Similarly, we may rewrite  the distribution 
$$S_\epsilon= \sum_{k=1}^3 \int_{\mathbb{R}^{3}} \partial_k\theta_\epsilon(y)  \big(u_k(x-y)-u_k(x)\big) \big( (\vec u(x-y)-\vec u(x))\cdot(\vec u_\epsilon(x)-\vec u(x)) \big)\, dy,$$ 
as
\begin{eqnarray*}
S_\epsilon&=& \sum_{k=1}^3 \int_{\mathbb{R}^{3}} \partial_k\theta_\epsilon(y)  
u_k(x-y) \vec u(x-y)\cdot\vec u_\epsilon(x)\, dy\\
&& -\sum_{k=1}^3 \int_{\mathbb{R}^{3}} \partial_k\theta_\epsilon(y)  
\left(u_k(x-y)-u_{k}(x)\right) \vec u(x-y)\cdot\vec u(x)\, dy\\
&&-\sum_{k=1}^3 \int_{\mathbb{R}^{3}} \partial_k\theta_\epsilon(y)  
u_k(x)\vec u(x-y)\cdot\vec u_\epsilon(x)\, dy\\
S_\epsilon&=& (\vec u\ast\theta_\epsilon)\cdot[\Div(\vec u\otimes\vec u)]\ast \theta_\epsilon -\sum_{k=1}^3 \int_{\mathbb{R}^{3}} \partial_k\theta_\epsilon(y)  
\left(u_k(x-y)-u_{k}(x)\right) \vec u(x-y)\cdot\vec u(x)\, dy\\
& &-\vec u_\epsilon\cdot[(\vec u\cdot\vec\nabla)\vec u_\epsilon].
\end{eqnarray*}
We thus have 
\begin{eqnarray*}
2S_\epsilon-R_\epsilon-\mu_{\varepsilon}&=&  -2\vec u_\epsilon\cdot[(\vec u\cdot\vec\nabla)\vec u_\epsilon]-\Div\left( \theta_\epsilon* (\vert\vec u\vert^2\vec u)-(\theta_\epsilon*\vert\vec u\vert^2)\vec u\right)+\Div(|\vu|^{2}\vu)\\
&=& -\Div(\vert\vec u_\epsilon\vert^2\vec u)-\Div\left(\theta_\epsilon* (\vert\vec u\vert^2\vec u)-(\theta_\epsilon*\vert\vec u\vert^2)\vec u\right)+\Div(|\vu|^{2}\vu).
\end{eqnarray*}
But, since $\vu\in L^3_t L^3_x(\Omega)$ we have 
$$ \lim_{\epsilon\rightarrow 0}  \theta_\epsilon* (\vert\vec u\vert^2\vec u)-(\theta_\epsilon*\vert\vec u\vert^2)\vec u=0\text{ in } L^1_{t}L^1_{x}(Q_{t_0,x_0,r/4}),$$
and $\Div(\vert\vec u_\epsilon\vert^2\vec u)\underset{\varepsilon\to 0}{\longrightarrow} \Div(|\vu|^{2}\vu)$ where the limit is taken in $\mathcal{D}'(Q_{t_0,x_0,r/4})$. Thus, we find that
$$ \lim_{\epsilon\rightarrow 0} 2S_\epsilon-R_\epsilon-\mu_\epsilon=0,$$ 
which is the formula given by Duchon and Robert. \hfill $\blacksquare$\\
%%%%%%%%%%%%%%%%%%%%%%%%%%%%%%%%%%%%%%%%%%%%%%

%%%%%%%%%%%%%%%%%%%%%%%%%%%%%%%%%%%%%%%%%%%%%%
%%%%%%%%%%%%%%%%%%%%%%%%%%%%%%%%%%%%%%%%%%%%%%

\quad\\[5mm]

\begin{flushright}
\begin{minipage}[r]{100mm}
Diego \textsc{Chamorro} (diego.chamorro@univ-evry.fr)\\
Pierre-Gilles \textsc{Lemari\'e-Rieusset} (plemarie@univ-evry.fr)\\
Kawther \textsc{Mayoufi} (kawther.mayoufi90@gmail.com)\\[3mm] 
Laboratoire de Math\'ematiques et Mod\'elisation d'\'Evry (LaMME), 
UEVE, UMR CNRS 8071\,\& ENSIIE\ 
\\[2mm]
 Universit\'e d'Evry Val d'Essonne\\
IBGBI, 23 Boulevard de France\\
91037 Evry Cedex, France
\end{minipage}
\end{flushright}

\end{document}